\documentstyle[12pt,amssymb,psfig,epsfig,color]{amsart}
\def\from#1\to{\mathpunct:#1\rightarrow}

\newtheorem{thm}{Theorem}[section]
\newtheorem{cor}[thm]{Corollary}
\newtheorem{lem}[thm]{Lemma}
\newtheorem{prop}[thm]{Proposition}
\theoremstyle{remark}

\newtheorem{example}[thm]{Example}

\theoremstyle{definition}

\numberwithin{equation}{section}
\numberwithin{figure}{section}

\def\note#1{}

\newcommand{\QED}{\rlap{$\sqcup$}$\sqcap$\smallskip}

\def\sss{\subsubsection}

\newcommand{\di}{\partial}

\newcommand{\ra}{\rightarrow}

\newcommand{\imply}{\Rightarrow}

\def\ssk{\smallskip}
\def\msk{\medskip}

\def\nin{\noindent}

\def\sm{\smallsetminus}

\newcommand{\cl}{\operatorname{cl}}
\newcommand{\inter}{\operatorname{int}}
\renewcommand{\mod}{\operatorname{mod}}

\newcommand{\tl}{\tilde}

\newcommand{\supp}{\operatorname{supp}}
\newcommand{\id}{\operatorname{id}}

\renewcommand{\d}{{\diamond}}

\def\Top{{\mathrm{Top}}}

\newcommand{\eps}{{\varepsilon}}

\newcommand{\de}{{\delta}}
\newcommand{\la}{{\lambda}}
\newcommand{\La}{{\Lambda}}
\newcommand{\si}{{\sigma}}

\newcommand{\Om}{{\Omega}}
\newcommand{\om}{{\omega}}

\newcommand{\ba}{{\mbox{\boldmath$\alpha$} }}
\newcommand{\sba}{{\mbox{\small\boldmath$\alpha$} }}
\newcommand{\bbe}{{\mbox{\boldmath$\beta$} }}

\newcommand{\AAA}{{\cal A}}
\newcommand{\BB}{{\cal B}}
\newcommand{\CC}{{\cal C}}

\newcommand{\FF}{{\cal F}}
\newcommand{\GG}{{\cal G}}

\newcommand{\HH}{{\cal H}}

\newcommand{\LL}{{\cal L}}

\newcommand{\NN}{{\cal N}}

\newcommand{\RR}{{\cal R}}
\newcommand{\SS}{{\cal S}}
\newcommand{\TT}{{\cal T}}

\newcommand{\WW}{{\cal W}}
\newcommand{\XX}{{\cal X}}
\newcommand{\YY}{{\cal Y}}
\newcommand{\ZZ}{{\cal Z}}

\newcommand{\D}{{\Bbb D}}

\renewcommand{\L}{{\Bbb L}}

\newcommand{\N}{{\Bbb N}}

\newcommand{\R}{{\Bbb R}}
\newcommand{\T}{{\Bbb T}}

\newcommand{\g}{{\bf g}}
\newcommand{\h}{{\bf h}}

\def\BS{{\mathbf{S}}}

\def\B0{{\mathbf{0}}}

\def\a{{\mathbf{a}}}
\def\b{{\mathbf{b}}}

\def\Crit{\mathrm{CV}}

\newcommand{\Comp}{\operatorname{Comp}}

\newcommand{\area}{\operatorname{area}}

\catcode`\@=12

\def\Empty{}
\newcommand\oplabel[1]{
  \def\OpArg{#1} \ifx \OpArg\Empty {} \else
  	\label{#1}
  \fi}
		
\newcommand{\comm}[1]{}
\newcommand{\comment}[1]{}

\def\makeabbrevs{%
\def\tA##1{$A_{##1}$}%
\def\ta##1{$a_{##1}^{-1}$}%
\def\tb##1{$b_{##1}^{-1}$}%
\def\o{\omega}\def\g{\gamma}\def\G{\Gamma}\def\h{\hat}\def\d{\delta}\def\D{\Delta}%
\def\O{\Omega}\def\b{\beta}\def\l{\lambda}\def\s{\sigma}\def\a{\alpha}}

\begin{document}

\bigskip\bigskip

\title[Quasi-Additivity Law]{The Quasi-Additivity Law \\ in Conformal Geometry }
\author {Jeremy Kahn and Mikhail Lyubich}
\address{Stony Brook University and University of Toronto} 
\date{\today}

\begin{abstract} 
On a Riemann surface $S$ of finite type containing a family of $N$ disjoint disks $D_i$ (``islands''),
we consider several natural conformal invariants measuring the distance from the islands to $\di S$ 
and separation between different islands. In a near degenerate situation we
establish a relation between them called the Quasi-Additivity Law.
We then  generalize it to a  Quasi-Invariance Law 
providing us with a transformation rule of the moduli in question under
covering maps.    
This rule (and in particular, its special case  called the  Covering Lemma) 
has important applications in holomorphic dynamics which will be addressed in the forthcoming notes. 
\end{abstract}

\comm{
We consider a Riemann surface $S$ of finite type containing a family of $N$ disjoint disks $D_i$,
and prove  the following Quasi-Additivity Law:
If the total extremal width $\sum \WW(S\sm D_i)$ is big enough (depending on $N$)
then it is comparable  with the extremal width $\WW (S,\,  \cup D_i)$
(under a certain ``separation assumption'') .  

We also consider a branched covering  $f: U\ra V$ of degree $N$ between two disks
that restricts to a map $\La\ra B$ of degree $d$ on some disk $\La \Subset U$.
We derive from the Quasi-Additivity Law that if $\mod(U\sm \La)$ is sufficiently small,  
then (under a ``collar assumption'')
the modulus is {\it quasi-invariant} under $f$, namely
$\mod(V\sm B)$ is comparable with  $d^2 \mod(U\sm \La)$.

This {\it Covering Lemma} has important consequences in holomorphic dynamics which will
be discussed in the forthcoming notes. 
}

\setcounter{tocdepth}{1}
 
\maketitle





\section{Introduction} 

Several central problems in holomorphic dynamics depend on the so-called {\it a priori} bounds,
that is,  uniform lower bounds on the conformal moduli of certain dynamically defined annuli.
So far, the only analytic tools suitable to this end (for unreal maps) were the basic properties of the
moduli of annuli (transformation rules and  the Gr\"otzcsh Inequality).
In this paper we design a new analytic tool, 
the {\it Covering Lemma},  that provides us, in a near degenerate situation,
with a much stronger 
version of the transformation rule for conformal moduli under covering maps. In the following papers, it is used  to generalize the
Yoccoz Theorem (on local connectivity of non-renormalizable Julia sets) to higher degree unicritical maps \cite{high} and to prove {\it a priori} bounds 
(and hence MLC) for
some classes of infinitely renormalizable quadratic maps \cite{K,decorations}. 
Further applications of this method (to multicritical maps) are under way, see \cite{KS,Chinese}.   

We will derive the Covering Lemma from a ``Quasi-Additivity Law'' relating three natural conformal moduli
for a Riemann surface with several Jordan disks marked. 
Let us formulate it precisely.

Let $S$ stand for a compact Riemann surface with boundary.
We denote the {\it extremal length}  of a family $\GG$ of curves by $\LL(\GG)$,
and we let $\WW(\GG)=\LL(\GG)^{-1}$ be the corresponding {\it extremal width} (see the Appendix).
Given a compact subset  $K\subset \inter S$, 
we let $\LL(S,K)$ and $\WW(S,K)$ be respectively the extremal length and
width of the family of curves in $S\sm K$ connecting $\di S$ to $K$.

An open subset $A\Subset  \inter S $  is called an {\it (open) archipelago}
if its closure is a Riemann surface of finite type (not necessarily connected)
with smooth boundary. Its connected components are called {\it islands}.
 
Let $A_j$ ($j = 1,\dots, N$) be a finite family  of archipelagos 
in $S$ with disjoint closures. 
We call the number 
$$ \Top= {\Top}_S \{A_j \} = -\chi(S) + \sum_j \#\Comp \partial A_j $$
the {\it topological  complexity} of the family of archipelagos.

 Let us introduce 3 conformal moduli of this family of archipelagos:
$$
  X =X_S\{A_j\}= \WW(S,\, \bigcup_{j=1}^N A_j);  
$$

\begin{equation}\label{Y}
   Y =Y_S\{A_j\} = \sum_{j=1}^N \WW(S, A_j),
\end{equation}

$$
   Z = Z_S\{A_j\}= \sum_{j=1}^N  \WW ( S \sm \bigcup_{k \not= j } A_k ,\, A_j ).  
$$
The first modulus measures the (inverse) extremal distance from the union of the archipelagos to the boundary of $S$,
the second one is the harmonic sum of the extremal  distances from the individual archipelagos to the boundary of $S$,
while the last one measures the (inverse) separation between the archipelagos. 

There are some obvious relations between these moduli: $X\leq Y\leq Z$ and $Y\leq N X$.
The goal of this paper is to establish one non-obvious relation in a  near degenerate situation, 
(i.e., when $Y$ is big), namely, to bound $Y$ by the geometric mean of $X$ and $Z$
with an absolute constant. The number $N$ of the archipelagos does not 
appear in the estimate: it only influences how degenerate the situation should be:

\proclaim Quasi-Additivity Law. 
  There exists $K$ depending only on the topological complexity of the
family of archipelagos such that: 
                    $$ Y\geq K \imply Y^2 \leq 2 XZ. 
\footnote{In fact, our proof shows that ``2'' can be replaced with any constant $C> 4/3$.
  On the other hand, one can show that $C< 32/27$ would not work.}
$$

The proof of this law will occupy most of the paper.

\msk\noindent{\bf A simple example.} 
\begin{figure}[htbp]
\begin{center}
\makeabbrevs
\input{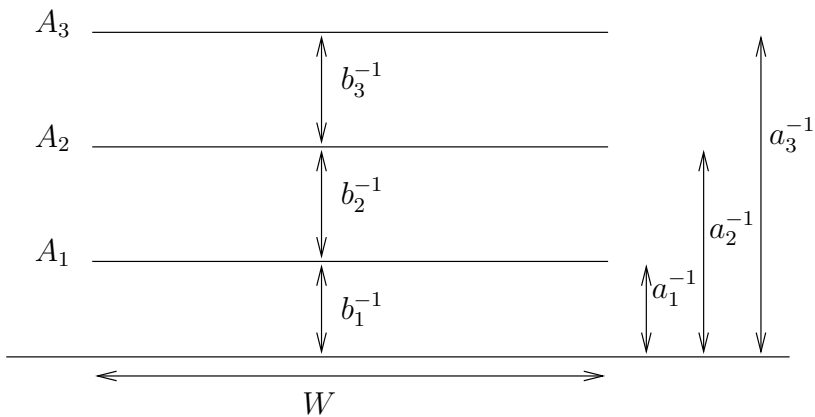}
\caption{Example: every island is a horizonal line segment.}
\label{example}
\end{center}
\end{figure}%
Figure \ref{example} presents a simple configuration of archipelagos
(consisting of a single island each) 
for which the asymptotics for the $X$-, $Y$- and $Z$-moduli 
can be calculated explicitly, so that the QA Law can be verified directly.
At the same time, this configuration is nearly optimal as the constant in the
QA Law is concerned.

Let $S$ be the closure of the upper half-plane in the Riemann sphere.
Given a sequence $a_1 > a_2 > \ldots > a_n > 0$,
let us consider archipelagos $A_i = [0, W] \times \{a_i^{-1}\}$, 
where $W$ is large in terms of the $a_i^{-1}$.
(Here our archipelagos are closed rather than open; 
see \S \ref{fractal} for a discussion).
Then
$$ X \sim W a_1,\quad Y \sim W \sum_{j=1}^n a_j,$$
and
$$
Z   \sim W \sum_{i=1}^n (b_i + b_{i+1}),
$$
where $b_i^{-1} = a_{i-1}^{-1} - a_i^{-1}$ (and $b_1 = a_1$, $b_{n+1} = 0$).
Then the QA Law in this case 
follows immediately from the arithmetic inequality
$$
\left(\sum_{j=1}^n a_j\right)^2 \le  \frac43b_1\sum_{j=1}^n b_j,
$$
which is proven in \S \ref{arithmetic sec}.

\ssk  Given $\xi\geq 1$, 
We say  that the archipelagos are $\xi$-{\it separated} if $Z\leq \xi Y$. 
The following immediate corollary shows that in a near degenerate situation,
under the separation assumption,   the moduli $X$ and $Y$ are comparable:

\proclaim QA Law with Separation.
Assume that the archipelagos $A_j\Subset \inter S$ 
are $\xi$-separated.
Then there exists $K$ depending only on $\xi$ and the topological complexity of the
family of archipelagos
 such that:
 $$ Y\geq K \imply Y\leq 2\xi X. $$


In \S \ref{variation} we  give several variations of the QA Law adapted to the
needs of  holomorphic dynamics.

\msk
We then generalize the QA Law to a Quasi-Invariance Law 
providing us with a transformation rule of the moduli in question under
covering maps in a near degenerate situation. Keeping in mind further applications,
we formulate in \S \ref{qi variations} a number of variations and special cases of this law.
Let us formulate here one of them. 

If we have a branched covering $f: U \ra V$ of degree $D$  between two disks
that restricts to a branched  covering  $f: \La\ra B$ of degree $d$ between smaller disks,
then a simple general estimate shows that $\mod(V\sm B) \leq D\mod(U\sm \La)$.
It turns out that given $d$, in a near degenerate situation 
the above moduli are, in  fact, comparable (under a ``collar assumption''):  

\proclaim Covering Lemma. 
Fix some  $\eta\in (0,1]$. 
 Let $U\supset \La'\supset \La$ and $V\supset B' \supset B$ be two nests of Jordan disks.
Let $f: (U, \La', \La) \ra (V, B', B)$ be a branched covering 
between the respective disks,  
 and let $ D  =  \deg(U\ra V)$, $d=\deg(\La'\ra B')$.
Under the following Collar Assumption:
$$ \mod(B'\sm B) > \eta \mod(U\sm \La),  $$ 
there exists an $\eps>0$ (depending on $\eta$ and $D$) such that such that
if  $$0< \mod(U\sm \La) < \eps$$ then 
$$ 
   \mod (V\sm B)  <  2\eta^{-1} d^2 \mod(U\sm \La).
$$

\begin{figure}[htbp]
\begin{center}
\makeabbrevs
\input{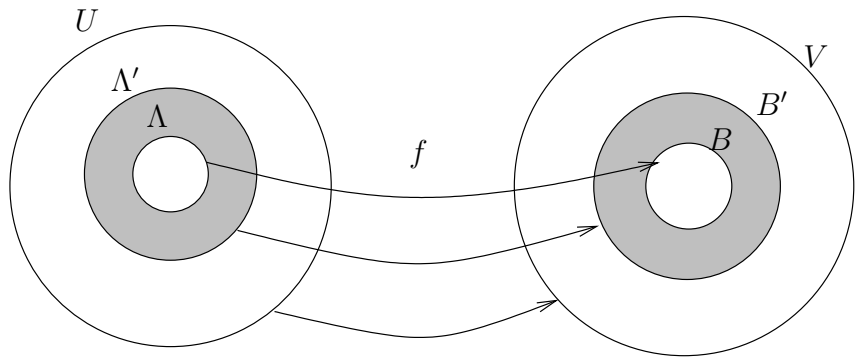}
\caption{Covering between two nests of three disks}
\label{three_disks}
\end{center}
\end{figure}

We derive the QI Law (and, in particular, this Covering Lemma) from the QA Law
by passing to an appropriate Galois covering of $U$.

The needed background in the extremal length techniques is summarized in 
the Appendix.

\msk {\bf Acknowledgment.}
 We thank Artur Avila and the referee for careful reading the manuscript and making many useful comments.
We  thank Nikita Selinger for giving us the best proof of Lemma \ref{harmonic sums}
(with the sharp constant).
We also thank all the Foundations that have supported this work: 
the Guggenheim Fellowship, Clay Mathematics Institute, NSF, and NSERC. 

\section{Quasi-Additivity Law}\label{quasi-add}


\subsection{Outline of the proof}

   Let us assume for simplicity (in this outline only) that $S$ and all the archipelagos
$A_j$  are  disks. Then each $S\sm A_j$ is an  annulus. Let us endow it with the vertical foliation
$\FF_j$ (the one that becomes genuinely vertical after uniformization of $S\sm A_j$ by a standard Euclidean
cylinder). Then $Y=\sum \WW(\FF_j)$. 

We begin  with analyzing topology of these foliations relative to our
family of archipelagos (\S\S \ref{paths and rectangles} - \ref{routes}). 
Namely, we associate to each leaf $\gamma$ of each $\FF_j$ a combinatorial invariant called its {\it route}. 
This invariant records the archipelagos visited by $\gamma$ (in order of first appearance) 
and some extra homotopy data
about $\gamma$. This data is selected in the minimal way 
to ensure that if two disjoint paths are   parallel (i.e., have the same route), then together
with appropriate arcs of the boundary of $ S\cup \bigcup A_j$ they bound a rectangle.
Moreover, any vertical path in this rectangle has the same route.
Thus, the vertical paths  with a given route vertically foliate a rectangle.
\note{``vertical'' is ambiguous} 

Let us consider one such rectangle, $P$, and let $(A_1,\dots, A_l)$ be the list of the archipelagos
visited by $P$.
This rectangle comes together with a sequence of 
associated ``big'' and ``little'' rectangles,  
$$
   P_k \subset S\sm \bigcup_{j= k}^N A_j, \quad Q_k \subset S\sm \bigcup_{j=1}^N A_j, \quad k=1,\dots, l.
$$
The big rectangles correspond to the pieces of its vertical boundary $\di^v P$ until its first entry to the archipelago
$A_k$, while the little ones correspond to the last piece of $\di^v P$ in $S\cup \bigcup A_j$.
The first of these rectangles, $Q_1$, is called ``initial''.
 
Cutting off from $P$ two buffers of width $4$ each, we obtain a truncated rectangle $\tl P$ coming
together with the associated truncated rectangles $\tl P_j$ and $\tl Q_j$.

At this point, we make use of a {\it Small Overlapping Principle} (\S \ref{so principle}) 
asserting  that families of curves with large extremal width have a relatively small intersection
(see \S \ref{so principle}). This implies that if two truncated little  rectangles overlap
(with matching vertical orientation) then the corresponding big rectangles have
comparable routes (i.e., one route is an extension of the other), see \S \ref{truncated rectangles}. 

This allows us to relate the moduli $X$ and $Z$ to the widths of the truncated small rectangles (\S \ref{forest}).
Namely, the total width of the truncated little rectangles is bounded by $2Z$, 
while the total width of the truncated initial little rectangles is bounded by $X$.  
On the other hand, the total width of the truncated big rectangles is bounded from below by
$(1-\de) Y$, as long as $Y > 16 s/\de$, where $s$ is the total number of the rectangles,
which can be bounded in terms of the topological complexity.

Moreover, by the Series Law for the extremal length, the width of each truncated big rectangle is bounded by the
harmonic sum of the widths of the associated little ones.  By an ``arithmetic inequality'' 
of \S \ref{arithmetic sec}, this yields the desired quadratic relation between the moduli $X$, $Y$ and $Z$.

\subsection{Paths and rectangles}\label{paths and rectangles} 
Let $S$ be a Riemann surface with  boundary.
All the curves  $\gamma: [0,1]\ra  S $ below will be considered naturally oriented. 
A curve $\gamma: [0,1]\ra  S $ is called {\it proper} if  $\gamma \{0,1\}\subset \di S$. 
Two proper curves are called {\it properly homotopic} in $S$  
if they are homotopic through a family of proper curves.
 A proper curve is called {\it trivial} if it is properly homotopic to a curve $[0,1]\ra \di S$. 
A {\it path}  in $S$ is a curve without self-intersections, 
i.e.,  an embedded (oriented) interval $[0,1]\ra S$.

In this paper, a {\it standard (Euclidean) rectangle} $E$ 
will mean $I\times [0,h]$ where $I$ is an interval of arbitrary type
(closed,  semi-closed, or open), and $h>0$. 
Its  {\it horizontal boundary}  $I\times \{0, h\}$ comprises the {\it base} $I\times \{0\}$ and
the {\it roof} $I\times \{ h\}$. 
A {\it vertical path} in $E$ is a path connecting its horizontal sides.
Every vertical path is naturally oriented (from the base to the roof) which endows $E$ with {\it vertical orientation}. 
The intervals  $\{x\}\times [0,h]$ will be referred to as {\it genuine vertical paths} in $E$;
together, they form the {\it genuine vertical foliation}.

A {\it (topological) rectangle} $P$ on a surface $S$ will mean an embedded Euclidean rectangle,
coming together with all the previously described affiliated structure: the horizontal boundary $\di^h P$ 
comprising the base and the roof, and the vertical orientation.
In what follows we will often deal with {\it properly embedded} rectangles, i.e., such that 
$\di^h P\subset \di S$. Any topological rectangle can be conformally uniformized by a standard
rectangle, supplying the former with the genuine vertical foliation. 

Similarly, a standard cylinder will mean $C=\T\times [0,h]$, 
where $\T$ is a round circle, coming together with the base and the roof, and the vertical orientation
(and the genuine vertical foliation, too). 
A {\it (topological) annulus} $R$ on $S$ is an embedded cylinder supplied with all the 
affiliated structure.

If we cut the annulus along  two disjoint vertical paths, we obtain two rectangles. 
This situation is special, as only one rectangle would  be cut off from any other Riemann surface: 

\begin{lem}\label{alpha-om}
 Assume $S$ is connected and not an annulus.
Let $\CC^1$ and $\CC^2$ be two disjoint properly homotopic non-trivial paths in $S$ such that $\inter \CC^i\subset \inter S$.
 
 \ssk\nin {\rm (i)}
Then there exist two unique arcs $\alpha$ and $\om$ on the boundary $\di S$ which together
with the paths $\CC^i$ bound a closed rectangle $P$ with base $\alpha$ and roof $\omega$. 

\ssk\nin {\rm (ii)}
 Let $(\CC^t)$, $1\leq t\leq 2$, be a proper homotopy between the above  paths,
   and let $(e^t)\subset \di S$ be the corresponding motion of the endpoint $e^t$ of $\CC^t$.
   Then the curve $(e^t)_{1\leq t\leq 2}$ is homotopic in $\di S$ rel the endpoints 
    to the arc $\om$ oriented from $e^1$ to $e^2$.      

\ssk\nin {\rm (iii)} Let $\CC^3$ be a third path which is disjoint and properly homotopic to the above two.
Let $P_j$, $j=1,2,3$, be the rectangles bounded by  the pairs of these three paths.
Then one of these rectangles is tiled by the other two. 

\end{lem}

\begin{pf} 
(i) 
 Let us consider the universal covering $\pi: \hat S\ra S$ of $S$.
It is conformally equivalent to $\Bar \D\sm K$, where $\bar \D$ is the closed unit disk and $K\subset \T$
is a nowhere dense compact subset of the unit circle (the limit set of the Fuchsian group of deck transformations).
Since the paths $\CC^i$ are properly homotopic, they lift to
(disjoint)  properly homotopic  paths $\hat \CC^i$ in $\hat S$.
Let these lifts begin at points $b^i\in \T$ and end at points $e^i\in \T$.
Then $b^1$ and $b^2$ (resp., $e^1$ and $e^2$) bound an arc $\hat\alpha \subset \di \hat S$ 
(resp. $\hat\om \subset \di \hat S$).
These two arcs are disjoint since the paths $\CC^i$ are non-trivial. 
They are also disjoint from the $\inter\CC^i\subset \inter \hat S$. 
Hence the four  paths, $\CC^1$, $\CC^2$, $\hat\alpha$ and $\hat\om$, bound  a closed rectangle $\hat P$ in $\hat S$.

Let us consider all the lifts $\hat\CC_j^i$ of $\CC^i$ that cross $\hat P$, where $\hat \CC^i_0\equiv \hat\CC^i$.
For each  $i=1,2$, the lifts $\hat C_j^i$ are pairwise disjoint 
since the paths $\CC^i$ do not have self-intersections.
Any two  paths $\hat \CC^1_j$ and $\hat \CC^2_k$ are disjoint as well
since $\CC^1$ and $\CC^2$ do not cross each other.
Hence each $\hat \CC^i_j$ is completely contained in $P$ and moreover, 
$\di \hat\CC^i_j\subset \hat\alpha\cup \hat\om$. But $\di\hat\CC^i_j$
cannot belong to one horizontal side, $\alpha$ or $\om$,
since the paths $\CC^i$ are non-trivial.
Thus, we obtain a  family of disjoint vertical paths $\hat\CC^i_j$ in $\hat P$. 

If one the above curves, say $\CC^1$, has  more than one lift, then let us consider the lift
$\hat\CC^1_1$ such that there are no other lifts in between $\hat\CC^1_0$ and $\hat \CC^1_0$. 
Then $\hat\CC^1_0$ and $\hat\CC^1_1$, together with two subarcs of $\hat\alpha$, 
and $\hat\om$ bound a rectangle $\hat \Pi$. The projection  of this rectangle to $S$
is a clopen annulus $R$ in $S$. Since $S$ is connected, $S=R$ contradicting our assumption.   

\ssk Thus, each  curve $\CC^i$ has only one lift to $\hat P$, so $\hat P\cap \pi^{-1}(\CC^i)=\hat \CC^i$. 
It follows that the paths $\CC^i$ lie on the boundary of $P \equiv  \pi(\hat P)$. 
Hence $\pi(\di\hat P) \subset \di P$,  
and the map $\pi: \hat P\ra P$ is proper.
Moreover, it is injective over  $\CC^i$ and hence has degree 1. 
Thus,  the map $\pi: \hat P\ra P$ is a homeomorphism.

\ssk
If there were two rectangles $P^1$ and $P^2$ as above then they would be glued along the paths $\CC^i$
to form an annulus.  

\ssk (ii)
  The homotopy $(\CC^t)$ lifts to a proper homotopy $\hat \CC^t$ on $\hat S$ between the lifts $\hat \CC^i$
considered in (i). The endpoint $\hat e^t$ of this lift moves along the component $\hat \xi$ of $\di \hat S$.
Since $\hat\xi$ is an interval, the curve $(\hat e^t)$ is homotopic to the arc $\hat \om$  on $\hat \xi$ rel the endpoints. 
Hence $(e^t)$ is homotopic to $\om$  on $\di S$ rel the endpoints. 

\ssk (iii)
 The paths $\CC^i$ lift to proper paths $\hat \CC^i$ in  $\hat S$ that begin and end on the same component of $\di \hat S$.
Then one of the lifted rectangles $\hat P_j$ is tiled by the other two.
Since $\pi: \hat P_j\ra P_j$ is a homeomorphism, the same is true for the $P_j$'s.

\comm{
\ssk (ii) 
Assume now that the paths $\CC^i$ are properly homotopic rel an arc $\om \subset \di S$ bounded by their endpoints.
Let $\beta$ be the component of $\di S$ containing $\om$. Let us double $S$ along $\beta$ and then cut the double
along $\om$.  We obtain a Riemann surface ${\bf S}= S\cup S'$ with $S\cap S'= \eta$,
where $\eta$ is the complement of $\om$ in  $\beta$
Moreover,  the paths $\CC^i$ are properly homotopic in $\bf S$. 
By (i), these paths together with appropriate arcs of $\di {\bf S}$ bound a rectangle $P\subset {\bf S}$.
Since the paths $\CC^i\subset S$ do not cross $\eta$, $\inter P\subset {\bf S}$ cannot intersect $\eta$ either,
and hence $P\subset S$. It follows that $\om$ (rather than the complementary arc $\om'\subset S'$)
 is a side of $P$.}   
\end{pf}

Somewhat loosely, we will say that the above rectangle $P$ is {\it bounded} by the curves $\CC^1$ and $\CC^2$.  

\comm{*** According to our convention, the sides $\alpha$ and $\om$ of  $P$ are
called  horizontal, while the sides $\CC^i$  are called  vertical.
Moreover, we assume that the paths $\CC^i$ begin on $\alpha$ and end on $\om$.  
Then  $\alpha$ is called the {\it base} of $P$, 
while $\om$ is called the {\it roof}.

Note that $P$  is endowed with a {\it vertical orientation}, i.e., the orientation of the vertical paths in $P$
consistent with the orientation of its vertical sides (namely, a positively oriented vertical path
begins on $\alpha$ and ends on  $\om$.) **** }

To avoid the ambiguity in the choice of the rectangle $P$,
{\it in what follows we assume that the Riemann surface $S$ under consideration is not an annulus.}
A simple trick shows that this assumption  does not reduce generality (see \S \ref{Vert fol}). 

Let us consider an archipelago $A$ in $S$.
Given a   proper path $\CC$ in $S$ that crosses  $\bar A$, 
let $a$ be the last point of intersection of $\CC$ with $\bar A$,  
and let $\de\subset S\sm A$ be the {\it terminal} closed segment of $\CC$  which connects $a$ to $\di S$.
Note that $\inter\de \subset \inter( S\sm A)$.  
If we have several paths $\CC^i$ as above,
we naturally label the corresponding objects as $a^i$ and $\de^i$, etc. 

\begin{figure}[htbp]
\begin{center}
\makeabbrevs
\input{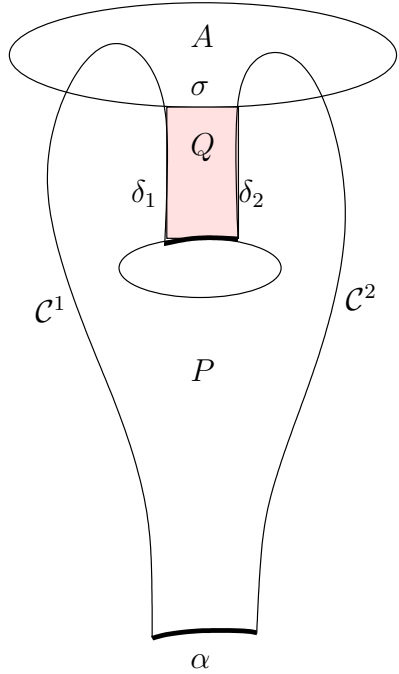}
\caption{Strange configuration of rectangles}
\label{small_wrong_picture}
\end{center}
\end{figure}

Two disjoint proper paths $\CC^1$ and $\CC^2$ in $S$ that cross $\bar A$  
are called {\it roof parallel} (rel $A$) if:
\begin{itemize}
\item  $\CC^1$ and $\CC^2$ are properly homotopic in $S$, and hence they bound a  ``big rectangle'' $P$; 
\item  
The  paths $\de^i$ are properly homotopic in $S\sm A$, and hence they bound a ``terminal little rectangle''
      $Q \subset S\sm A$;
\item 
    The rectangles $P$ and $Q$ share the roof (Figure \ref{small_wrong_picture} illustrates that this is not automatic.)
\end{itemize}

Two paths are called {\it base parallel} (rel $A$) if after reversing orientation
they become roof parallel. Initial segments of these paths bound an initial little rectangle $Q_1\subset S\sm A$
which shares the base with $P$. Two paths a called {\it parallel} if they are roof and base parallel.

We will now formulate several statements about roof parallel paths. 
The corresponding statements about base parallel paths are obtained by reversing orientation,
and the corresponding statements about parallel paths immediately follow.

\begin{lem}\label{intermediate paths}
Let $\CC^1$ and $\CC^2$ be two roof parallel (rel $A$) proper paths in $S$,
and $P$ and $Q$ be the corresponding big and little rectangles.
  Let $\CC$ be a positively oriented vertical path in $P$
which is disjoint from the $\CC^i$.
Then it is roof parallel (rel $A$) to each $\CC^i$. 
Moreover, its terminal segment  $\de$ is a vertical path in $Q$.
\end{lem}

\begin{pf}
 Any vertical  path in $P$ is properly homotopic to the sides $\CC^i$. 
Let $P^i$ be the big rectangles bounded by the paths $\CC$ and $\CC^i$, 
and let $\om^i$ be their roofs, $i=1,2$. 
Of course, they tile the roof $\om$, overlapping at the endpoint $e$ of $\CC$.

Let $\CC'$ be the path $\CC$ with reverse orientation.
Since $P$ and $Q$ share the roof,  some initial segment of $\CC'$  is contained in $Q$. 
Since $\CC'$ is proper, it must exit $Q$. Since $\inter \CC'$ is disjoint from the vertical sides and
the roof of $Q$, it can exit $Q$ only through its base, $\si$. Let $a$ be the first point of intersection 
between $\CC'$ and $\si$. Then the terminal segment $\de$ of $\CC$ that begins at $a$ 
is a positively oriented vertical path in $Q$. 
Hence it is properly homotopic in $S\sm A$ to the paths $\de^i$.

Let $Q^i\subset S\sm \bar A$ be the little rectangles bounded by the paths $\de$ and $\de^i$, $i=1,2$.
Since $\de$ is a vertical path in $Q$ ending at $e$, the arcs $\om_i$ are the roofs of the little rectangles $Q^i$.
Thus, $Q_i$ respectively share the roofs with $P_i$. 
\end{pf}

  The following lemma will be used for counting the number of parallel classes
(see Lemmas \ref{number of par classes} and \ref{route count}): 

\begin{lem}\label{3 paths}
  Let $\CC^i$ be three disjoint properly homotopic  paths in $S$ crossing the archipelago $\bar A$
in such a way that their terminal segments $\de^i$ 
are properly homotopic in $S\sm A$.
Then at least two of these paths are roof parallel rel $\bar A$. 
\end{lem}

\begin{pf}

 For $i=1,2,3$,  let $P_i$ be the big rectangle  bounded by the paths $\CC^k$ and $\CC^l$ with $\{i,k,l\}=\{1,2,3\}$,
and let $Q_i$ be the corresponding little rectangles. 
Let $\om_i$ be the roofs of the $P_i$, and let $\la_i$ be the roofs of the $Q_i$. 
We need to show that one of the roofs $\om_i$ coincides with the corresponding $\la_i$. 

Since by Lemma \ref{alpha-om} (iii) one of the big rectangles, say $P_1$, is tiled by the other two,
the roof $\om_1$ is tiled by $\om_2$ and $\om_3$. 
Denote the complements of the roofs $\om_i$ by $\om_i'$. 
If $\om_i\not=\la_i$ for $i=2,3$,
then $\la_2=\om_2'=\om_1'\cup \om_3$ and similarly $\la_3=\om_1'\cup \om_2$.
Hence $\la_2\cup \la_3=\om_1'\cup \om_2\cup \om_3=\eta$, where $\eta$ is the whole component of $\di S$
containing the endpoints of the paths $\CC^i$. But it is impossible since one of the roofs $\la_i$
is tiled by the other two (as one of the little rectangles $Q_i$ is tiled by the other two).  
\end{pf}

\begin{cor}\label{3 paths cor}
  Let $\CC^i$ be five disjoint properly homotopic  paths in $S$ crossing the archipelago $\bar A$
in such a way that their terminal and initial segments 
are (respectively) properly homotopic in $S\sm A$.
Then at least two of these paths are parallel rel $\bar A$. 
\end{cor}

Let us now enlarge the notion of parallel to an equivalence relation on the class $\AAA$ of all proper
curves $\CC$ in $S$ crossing the archipelago $\bar A$.
Let us say that two curves $\CC^1$ and $\CC^2$ of class $\AAA$ are {\it roof equivalent}
if 
\begin{itemize} 
  \item They are properly homotopic in $\CC$;
  \item The terminal segments $\de^1$ and $\de^2$ are properly homotopic in $S\sm A$;
  \item The motions of the endpoints,  $(e^t)$ and  $(q^t)$, of the above homotopies 
   are  homotopic  (rel endpoints) curves on $\di S$.
\end{itemize}

The definitions of {\it base equivalent} and {\it equivalent} paths are straightforward.
Again, we restrict ourselves to a statement concerning roof equivalence only:

\begin{lem}
     Two disjoint curves $\CC^1$ and $\CC^2$ of class $\AAA$ are roof parallel if and only if they are roof equivalent.  
\end{lem}

\begin{pf}
   If $\CC^1$ and $\CC^2$ are roof parallel then they are homotopic within the big rectangle $P$ in such a way
that the endpoint $e^t$ parametrizes the roof $\om$. Similarly, the curves $\de^1$ and $\de^2$ are homotopic 
in $Q$ in such a way that $q^t$ parametrizes the same roof $\om$. So, the motions of the endpoints are homotopic.

Vice versa, by Lemma \ref{alpha-om} (ii), the homotopy class of the endpoint motion  determines the roof of the
rectangle.  
 \end{pf}

{\it In what follows, (roof/base) equivalent curves (not necessarily disjoint) will also be called (roof/base) parallel.} 
Also, ``parallel in $S$ (rel~$\emptyset$)'' just means  ``properly homotopic'' in $S$.

\msk
We close with two combinatorial lemmas. 
\begin{lem}
  Suppose that $S$ is a Riemann surface of finite topological type such that each connected component of $S$
has negative Euler characteristic. 
  Then there can be at most $-3\chi(S)$ disjoint non-parallel (rel~$\emptyset$) proper paths in $S$
\end{lem}
\begin{pf}
  Removing boundary from $S$, we obtain a Riemann surface homeomorphic to a
a compact Riemann surface $\BS$ with finitely many punctures $v_k$, $k=1,\dots, n$,  
where say, the first $l$ of them  correspond to the removed boundary. 
Proper paths in $S$ correspond to  paths in $\BS\sm \{v_k\}$ connecting 
two of the first $l$ punctures.
Of course, if we allow ourselves to connect other vertices as well, we obtain only more paths.
So, we can assume in the  first place that $S=\BS\sm \{v_k\}_{k=1}^n$ and $l=n$
(and of course, we can assume that $n\geq 1$). 
Since the Euler characteristic is additive, we can also assume that $S$ is connected.

Let us call the punctures ``vertices'' and non-trivial paths in $\BS\sm \{v_k\}$ connecting them ``edges''.
It is well-known that any finite family $\FF$ of  disjoint non-parallel edges can be completed
to  a triangulation of the surface $\BS$ with the same vertices $v_k$ (provided $\chi(S)<0$). 
[To see it, let us first complete $\FF$ to a connected graph containing all the vertices $v_k$.  
Let us then consider any ``face'' $D$ of it, i.e., a component of the complement of the edges.
If $D$ has positive genus, we can add to $\FF$ a closed non-dividing edge connecting some vertex to 
itself. Cutting along this edge, we reduce the genus of $D$. Proceeding in this way,  
we will eventually obtain a graph whose faces are polygons.
None of these faces can be a bigon, since the edges are not parallel.
It cannot be a one- or zero-gon either since $\chi(S)<0$.  
Thus, all the polygons are at least $m$-gons with $m\geq 3$, 
and  we can further triangulate them.]
      
   Let us apply the Euler formula to this triangulation:
 $$
    F - E + V = \chi(\BS),
 $$
  where $E$ is the number of proper paths, $V = n$, and $3F = 2E$.
  Therefore $- E/3  = \chi(\BS)-n = \chi(S)$,
  and we are done.
\end{pf}

\begin{lem}\label{number of par classes}
  Suppose that $A$ is an archipelago on $S$,
  and let $\FF$ be a set of disjoint proper paths on $S$.
  Then there are at most 
  $$
  - 108\, \chi(S)\, \chi(S\sm A)^2
  $$
  distinct parallel classes (rel A) in $\FF$.
\end{lem}
\begin{pf}
  There are at most $-3\chi(S)$ distinct homotopy classes of curves $\gamma$ in $\FF$,
  and at most $-3\chi(S \sm A)$ distinct homotopy class 
  for the initial and final segments of $\gamma$.
  By  Corollary \ref{3 paths cor},  there are at most four distinct parallel classes,
  given the homotopy classes for $\gamma$ and its initial and terminal segments.
\end{pf}

\subsection{Routes and associated rectangles}\label{routes}
  Let us now consider  a finite family $\AAA$ of archipelagos $A_j$ ($j=1,\dots, N$)
in $S$ with disjoint closures. 
Let us consider a path $\CC$  in $S$ that begins at $b\subset \di S$
and ends at  a point $e$ on some archipelago $\bar A $. 
Such a path is called {\it good} if $\inter \CC$ does not intersect $\di S\cup \bar A$.

Given a good path $\CC$ in $S$, let us relabel (if needed) our archipelagos so that
$( A_1, \dots,  A_l\equiv A )$ is the sequence of distinct archipelagos whose closures are crossed by $\CC$
ordered according to their first appearance, while  $A_{l+1}, \dots, A_N$ are the archipelagos
that are not crossed by $\CC$ ordered in an arbitrary way. 
Thus, for any $1\leq i<j \leq l$,  the  path $\CC$ enters $A_i$ for the first time before it enters $A_j$.  
Note that though $\CC$ can enter each archipelagos $A_i$ ($1\leq i\leq l$)
many times, it is recorded only once. 

Let $e_j$ be the first point of intersection of $\CC$ with $\bar A_j$, and let 
$\CC_j$ be the segment of $\CC$ bounded by $b\equiv e_0$ and $e_j$.  
In this way we obtain the {\it associated sequence}
$$\CC_1\subset \CC_2\subset \dots  \CC_l\equiv \CC $$ 
of good paths in $S$.
We let   $|\CC| =l$ and call it   the {\it height}  of  $\CC$.

\begin{figure}[htbp]
\begin{center}
\makeabbrevs
\input{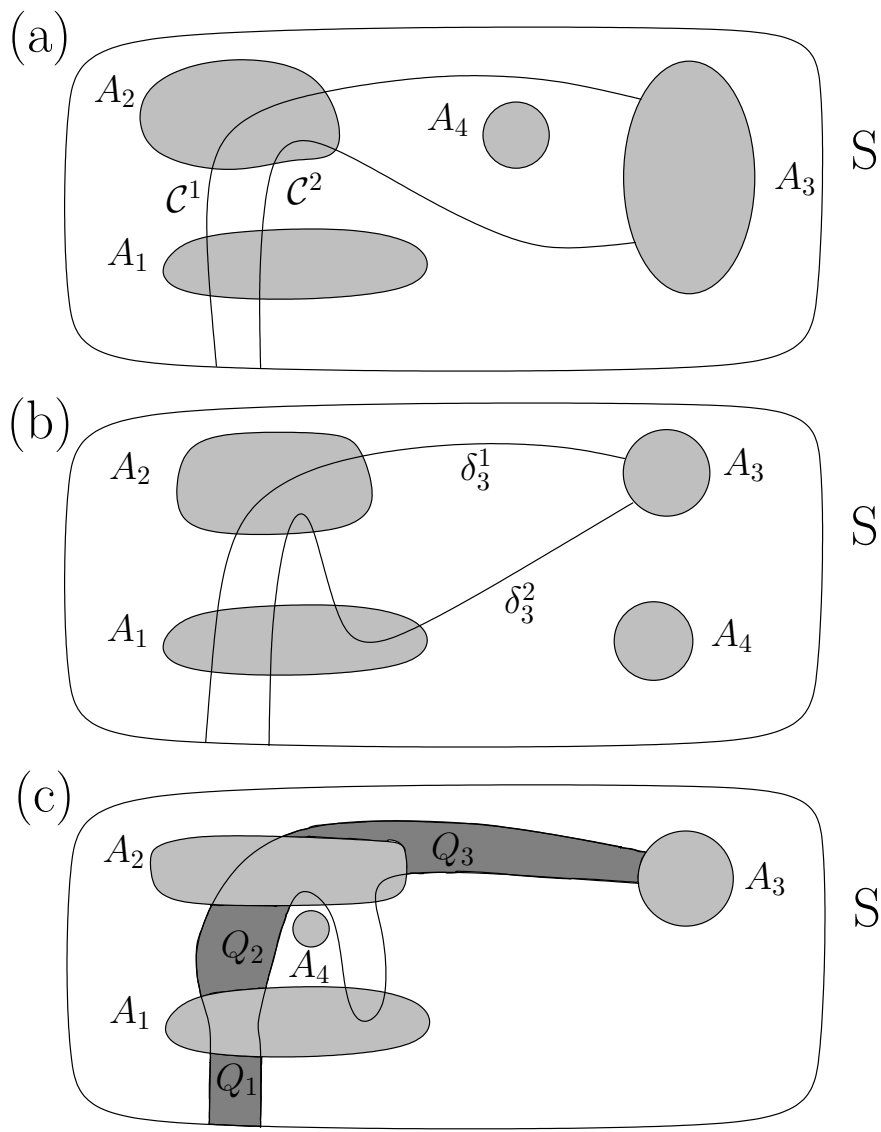}
\caption{
This picture illustrates the notion of parallelism.
Here the family $\AAA$ comprises four archipelagos $A_i$ each consisting of a single island.
The routes of the paths $\CC^1$ and $\CC^2$ have height $l=3$, and
$\La_3= A_3\cup A_4$, $\Om_3= A_1\cup A_2$.
The paths on figure (a) are not parallel since they are not properly homotopic in $S\sm \La_3$.
The paths on figure (b) are not parallel since their terminal arcs, $\de^1_3$ and $\de^2_3$,
are not properly homotopic in $S\sm \La$.
On the other hand, the paths on (c) are parallel, notwithstanding $\CC^2$ visits the island $A_1$ twice,
while $\CC^1$ visits it only once. 
In all three cases, the initial segments of the paths (of height two), $\CC^1_2$ and $\CC^2_2$,
are obviously parallel. 
}
\label{three_cases}
\end{center}
\end{figure}

Let 
$$ 
       \La_j= \bigcup_{i=j}^N A_i, \quad   \Om_j = \bigcup_{i=1}^{j-1} A_i. 
$$ 
(Note that $\Om_1=\emptyset$. Also, we let $\La\equiv \La_1$ be the union of all archipelagos.)
Then $\CC_j$ is a proper path in $S\sm \La_j$, 
and $\Om_j$  is an archipelago in $ S\sm \La_j$.   
Let $\ba_j$ be the class of proper paths in $S\sm \La_j$ parallel to $\CC_j$ rel $\Om_j$. 
We say that these paths and classes are  {\it associated} to $\CC$.
The sequence of the associated parallel classes,
$$
  \RR(\CC)= (\ba_j)_{j=1}^l,
$$
is called the {\it route} of $\CC$. Note that the route determines the base component  of $\di S$ 
where $\CC$ begins,  
and  the components of  $\di A_j$ where the curves $\CC_j$ end.
Two good paths are called {\it parallel} rel the family $\AAA$ of archipelagos  if they have the same route.
Note that parallel paths can cross some particular archipelagos $A$ different number of times
(see Figure \ref{three_cases}). 

\msk
We will now derive a bound on the number of routes: 
 
\begin{lem}\label{route count} 
Let $A_1, \ldots, A_N$ be distinct archipelagoes in $S$,
and let $\TT$ be a set of disjoint good paths in $S$,.
Then among the elements of $\TT$ there are at most $s(Top, N) = N!\, (108\Top^3)^{N+1}$ 
distinct routes rel $\{ A_j\}$
(where $\Top = \Top_S\{A_j\}$ is the topological complexity defined in the introduction).
\end{lem}
\begin{pf}
{
  \def\NN{\frac{N!}{(N-k)!}}
  Let us bound the number of routes $\RR(\CC)$ (for $\CC \in \TT$) 
  for which $A_1, \ldots A_k$ are visited in sequence
  (so that the terminal point of $\CC_j$ lies in $\bar A_j$).
  By the previous Lemma, 
  there are at most 
  $$
  108\, \chi(S \sm \La_j)\, \chi(S \sm \La)^2  \le 108\Top^3
  $$
  distinct parallel classes for $C_j$,
  so there are at most 
  $(108\, \Top^3)^k$
  distinct routes which visit $A_1, \ldots, A_k$ in sequence.
  There are $\NN$ injective functions 
  $\sigma: \{1, \ldots, k\} \to \{1, \ldots, N\}$, 
  so there are at most 
  $$\NN(108\Top^3)^k$$
  distinct routes of length $k$. 
The total number of these routes is bounded as desired.
}
\end{pf}

Let us consider two disjoint parallel good paths $\CC^1$ and $\CC^2$ with route of height $l$.
By Lemma \ref{alpha-om}, these two paths, together with a base  $\alpha$ and a roof $\om$,
bound a {\it good big rectangle} $P$.    
Moreover, for each $j=1,\dots, l$,
the associated good  paths $\CC^1_j$ and $\CC^2_j$, together with a base  $\alpha_j$ and a roof $\om_j$,
bound an {\it associated good big rectangle} $P_j\subset S\sm \La_j$, where $P_l \equiv P$.
In fact, the $P_j$ share the same base, i.e. $\alpha = \alpha_j$, since they share the base with the same 
{\it associated initial little rectangle} $Q_1 \equiv P_1$.  
Furthermore, each rectangle $P_j$ shares the roof with 
{\it associated (terminal) little rectangle} $Q_j$, $j=2,\dots, l$, 
bounded by the terminal paths $\de_j^1$ and $\de^2_j$, a base $\si_j$,  and the roof $\om_j$.
Note that the little rectangles $Q_j$ are not necessarily contained in the  big rectangle $P$ (see Figure \ref{big_little}).
All the above rectangles are vertically orientated.

\begin{figure}[htbp]
\begin{center}
\makeabbrevs
\input{figures/big_little2.pstex_t}
\caption{}
\label{big_little}
\end{center}
\end{figure}VC
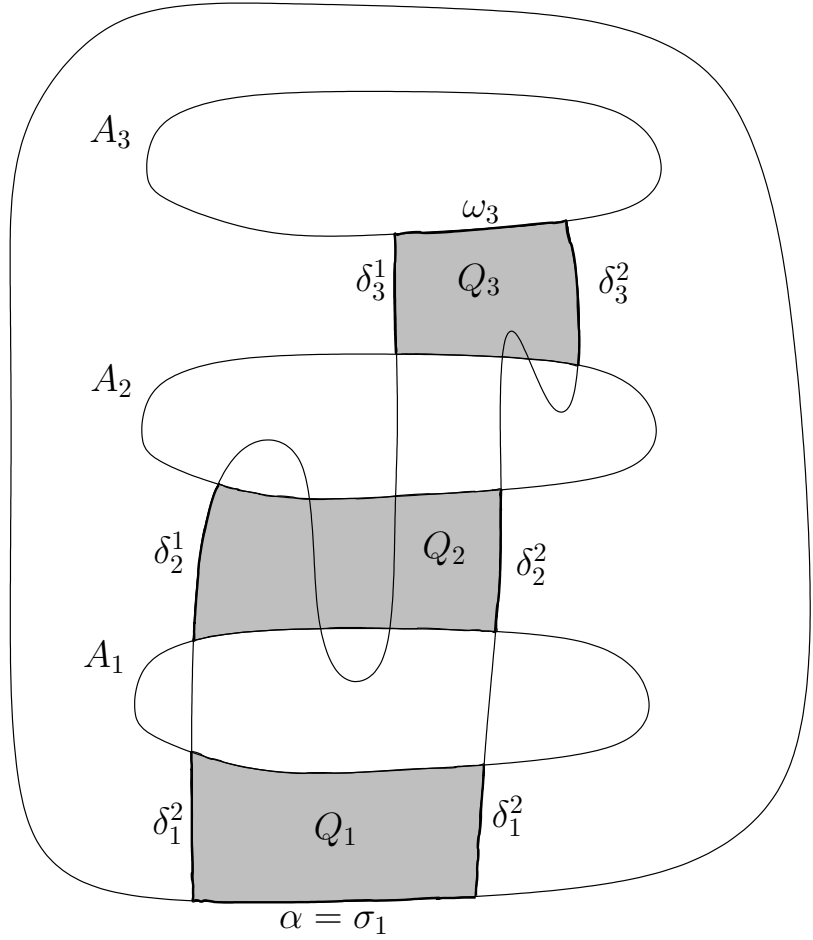
{This picture illustrates that the little rectangles $Q_i$ (shaded)
   are not necessarily contained in the big ones.}

We say that a path $\CC$ {\it (positively) vertically overflows} a little rectangle $Q_j$ if $\CC$ contains a segment $\de$
which is a (positively oriented) vertical path in $Q_j$.


\begin{lem}\label{good rectangles}
Let $\CC^1$ and $\CC^2$ be two disjoint parallel (rel $\AAA$) good paths of height $l$, 
and let $P\equiv P_l$ be the corresponding  good big rectangle.
Let $\CC$ be a positively oriented vertical path in $P$. 
Then it is parallel to $\CC^1$ and $\CC^2$ (rel $\AAA$)  
and, in particular, it has height $l$. 
Moreover,  $\CC$ positively vertically overflows all associated little rectangles $Q_j$, $j=1,\dots, l$. 
\end{lem}

\begin{pf}
Let us begin with the last assertion. 
For $j=l$ and $j=1$ it immediately follows from Lemma \ref{intermediate paths}
(by reversing orientation for $j=1$).
Let $1<j <l$. 
 Since $P_j$ has the same base $\alpha\subset \di S$ as $P$,
a little initial segment of $\CC$ is contained in $P_j$.
On the other hand, the endpoint of $\CC$ belongs to the archipelago $\bar A_l$
which is disjoint from  $P_j$ since
$$ 
     P_j\subset S\sm \La_j \subset S\sm \bar \La_l.
 $$  
Hence the curve $\CC$ must exit the rectangle $P_j$.
But since  $\CC$  is a vertical curve in $P$, it can exit $P_j$ only through the roof $\om_j$. 
Let $e_j$ be the first intersection point of $\CC$ with this roof. 
Then the initial segment $\CC_j$ of $\CC$ with endpoint $e_j$ is a vertical path of $P_j$.    
By Lemma \ref{intermediate paths}, it positively vertically overflows the little rectangle $Q_j$.
All the more, $\CC$ does so. 

Since each $P_j$ is a good big rectangle as well, we can apply to it the previous result
and conclude that for any $i\leq j$, $\CC_j$ vertically overflows $Q_i$.
In particular it crosses the roof $\om_i\subset \di A_i$, and hence $\CC_i\subset \CC_j$.  

Let us show that $\CC_1\subset\dots \subset \CC_l$ is the associated sequence of good paths. 
Since all the paths $\CC_j$ are good initial segments of $\CC$, 
it is  part of the associated sequence.
Moreover, $\CC$ does not contain any other good initial segment since all other 
archipelagos $A_k$, $k=l+1, \dots, N$, are disjoint from $P$.  

In particular, $\CC$ has the same height $l$ as $\CC^1$.
Moreover, by Lemma~\ref{intermediate paths}, the paths  $\CC_j$ are parallel to $\CC^1_j$ and $\CC^2_j$
rel $\Om_j$. Hence $\CC$ is parallel to $\CC^1$ and $\CC^2$ rel $\AAA$. 
\end{pf}

The previous lemma can be sharpened as follows: 

\begin{lem}\label{good rectangles-2}
Let $\CC^1$ and $\CC^2$ be two disjoint parallel (rel $\AAA$) good paths of height $l$, 
and let $P\equiv P_l$ be the corresponding  good big rectangle with base $\alpha$.
Let $\CC$ be a good path disjoint from $\CC^1$ and $\CC^2$ which begins on $\alpha$. 
Then  either the  route $\RR(\CC)$  extends $\RR(\CC^1)=\RR(\CC^2)$, or the other way around. 
\end{lem}

\begin{pf}
  Assume $\CC$ is not contained in the rectangle $P$. Then it must exit $P$ through the roof $\om$.  
Let $e$ be the first point of intersection of $\CC$ with $\om$. Then the initial   segment $\CC^*$  of $\CC$ ending at $e$
is a vertical path in $P$. By Lemma \ref{good rectangles},  $\RR(\CC^*)= \RR(\CC^1)$, 
so that $\RR(\CC)$ extends $\RR(\CC^1)$. 

\ssk
Assume now that $\CC\subset P$.
  Let us consider the biggest $j\leq l$ such that $\CC$ intersects  the roof $\om_j$ of the good big rectangle $P_j$,
and let $e_j\in \CC\cap \om_j$ be the first intersection point.  Then the initial segment $\CC_j$ of $\CC$
with endpoint $e_j$ is a vertical path in $P_j$. By Lemma \ref{good rectangles},
it has the same  route as $\CC_j^1$. In particular, it crosses all the archipelagos $A_i$, $i=1,\dots, j$. 

But in fact, $\CC=\CC_j$, for otherwise $\CC$ (being good) would end at some archipelago $A_i$ with $i>j$.
For $i>l$ it is impossible since those archipelagos are disjoint from $P$. For $i\in [j+1, l]$ it is impossible
for otherwise $\CC$ would exit the rectangle $\inter P_i$ and hence would cross the roof $\om_i$.   

We conclude that  $\RR(\CC^1)$ is an extension of $\RR(\CC_j)=\RR(\CC)$. 
\end{pf}

Let us now consider two disjoint vertical curves $\Gamma^1$ and $\Gamma^2$ in a good rectangle $P$.
Together with appropriate  base and roof arcs, they bound a truncated good  rectangle $\tl P\subset P$.

\begin{lem}\label{inclusions}
 For the associated sequence of  little rectangles, we have:
  $\tl Q_j \subset Q_j $. 
\end{lem}

\begin{pf}
By   Lemma \ref{good rectangles}, $\Gamma^1$ and $\Gamma^2$ have the same route as $P$.
Let us consider the associated sequences of good curves $\Gamma^1_j$ and $\Gamma^2_j$, $j=1,\dots, l$.
Let $\tl\de^1_j$ and $\tl\de^2_j$ be the terminal paths in $S\sm \bigcup A_j$  of these curves.
By definition, $\tl Q_j$ is the rectangle  bounded by these two paths, 
together with two appropriate horizontal arcs. 
By   Lemma \ref{good rectangles}, the $\tl\de^i_j$ are vertical paths in  
the little rectangle $Q_j$.
Hence $\tl Q_j\subset Q_j$.
\end{pf} 

Finally, we have the following important disjointness property: 

\begin{prop}\label{disjointness}
   Let $P$ and $P'$ be two good rectangles 
with disjoint vertical boundaries. 
Assume that some associated little rectangles, $Q_j$ and $Q_k'$, have
a non-trivial overlap. 
Then they represent the same proper homotopy class in $S\sm \La$ (up to orientation).
If their orientations match,
then one of the routes, $\RR(P)$ or $\RR(P')$, is an extension of the other,
and $j=k$.
\end{prop} 

\begin{pf}
Since the overlapping little rectangles $Q_j$ and $Q_k'$ have  disjoint vertical  boundaries, 
one of the vertical  boundary components, say  $\de_k' \subset \di Q_k'$, 
must be a vertical path in the other rectangle, $Q_j$, which implies the first assertion. 

Assume the vertical orientation of $Q_j$ and $Q_k'$ match.
Let $\CC'$ be the vertical boundary component of $P'$ containing the path $\de_k'$,
and let $\CC_k'$ be the associated good curve ending with the path $\de_k'$.

Let us consider the associated with $P$ good big rectangle $P_j$ 
(with the little rectangle  $Q_j$ just under its roof $\om_j$).
Since the path $\de_k'$ is positively oriented in $Q_j$,
it ends on the roof $\om_j$. Thus, the whole curve $\CC_k'$ also ends on $\om_j$. 
But  since $\CC_k'$ is good, its interior does not cross $\om_j$.
Neither can it cross the vertical boundary of $P_j$ (by the assumption).
Hence $\CC_k'$ is trapped in $P_j$, and must begin on the base $\alpha_j$ of $P_j$. 

Thus, $\CC_k'$ is a vertical curve in $P_j$.
By Lemma \ref{good rectangles}, $\CC_k'$ and $P_j$ have the same height, so that $k=j$.  
By Lemma \ref{good rectangles-2}, the route $\RR(\CC')=\RR(P')$ is either an extension of $\RR(P)$,
or the other way around.  
\end{pf}

\subsection{Harmonic foliations}\label{Vert fol} 

Let now $\BS$ be a compact Riemann surface with boundary, and let $S$ be obtained from $\BS$ by
making finitely many punctures $p_k\in \inter \BS$. We let $\di S= \di \BS$.  


By making a few artificial punctures
(depending only on the topological complexity of the family of archipelagos),
we can  ensure that {\it no component of $S\sm A_j$ is an annulus} 
(see our convention after~Lemma \ref{alpha-om} and Figure \ref{long_island}).
Note that making extra  punctures  does  not change extremal lengths of the path
families in question.

\begin{figure}[htbp]
\begin{center}
\makeabbrevs
\input{figures/long_island2.pstex_t}
\caption{}
\label{long_island}
\end{center}
\end{figure}VC
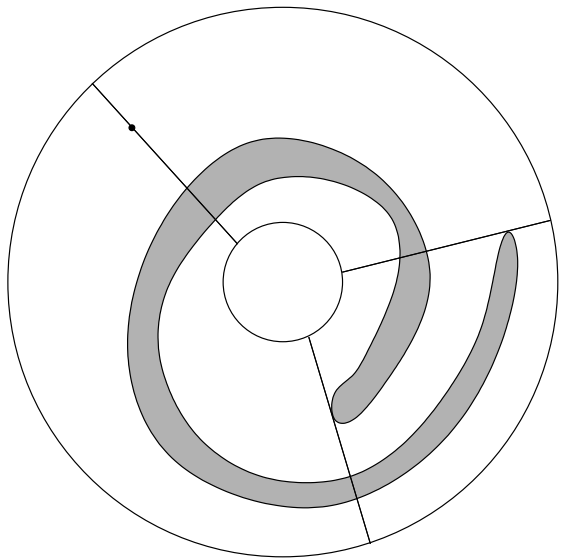
{Long Island. On this picture, $S$ is an annulus with one island on it. 
Without an artificial puncture,
all the leaves of the harmonic foliation would be in the same parallel class. 
With the puncture, 
the leaves are decomposed into three parallel classes that form three rectangles.}

 Let us consider the harmonic measure $\om_j(z)= \om_{S\sm A_j} (\di A_j, z)$ of  $\di A_j$ in the 
Riemann surface $\BS\sm A_j $ (see \cite{A}). It is a unique harmonic function on $\inter (\BS\sm A_j)$ 
equal to 1 on $\di A_j$ and vanishing on $\di \BS$. 
For instance, if $\BS$ and $A_j$ are disks, then $\om_j$ is the  height function on the annulus $\BS\sm A_j$ 
uniformized  by the flat cylinder $C_j$ with height 1 in such a way that $\di \BS$ is the base of it. 

The {\it harmonic foliation} $\FF_j$ on $\BS$ is the phase portrait of the gradient flow $\gamma_j^t$ of $\om_j$.
It has finitely many saddle type singularities
(with finitely many incoming and outgoing separatricies),
where the punctures are considered to be  singularities as well. 
It is  oriented according to the direction of the gradient flow.
Each non-singular leaf of $\FF_j$ begins on $\di \BS$ and ends on $\di A_j$. 
In the case when $\BS$ is a topological annulus, 
$\FF_j$ is the genuinely vertical foliation on the uniformizing cylinder $C_j$. 

Let us remove  from $S\sm A_j$ all separatricies $O^k $ of the foliation $\FF_j$
and take the the components of  $S\sm (A_j\cup\bigcup O^k)$.
We obtain finitely many  rectangles $\Pi=\Pi_j^m$ 
foliated by the harmonic leaves. 
Indeed, take some component $\la$ of $\di S\sm \bigcup O^k$.
The gradient flow brings every point $z\in \la$ in time 1 to some archipelago $ A_j$,
and these trajectories fill in some component $\Pi$ of $S\sm A_j\sm \bigcup O^k$.
The map 
$$ 
      (z, t)\ra (z, \gamma_j^t (z)), \quad z\in \la, \ t\in [0,1],
$$ 
provides us with the rectangular structure on $\Pi$.  
(Since every annuli component of $\BS\sm A_j$ contains a puncture, there are no annuli among the $\Pi_i$'s.)

The conjugate harmonic function $\om_j^*$  induces the natural  transverse measure
 on the  $\Pi_j^m$. 
In fact,   
the map $ \om_j+ i \om_j^*$
 provides us with the  uniformization of $\Pi^m_j$  by a standard rectangle of height 1.  

Every rectangle $\Pi_j^m$ represents some non-trivial proper homotopy class of paths in $S\sm A_j$. 
Moreover, different rectangles represent different classes.
Indeed, if two leaves, $\gamma$ and $\gamma'$, of $\FF_j$ are properly homotopic in  $S\sm A_j$,
then  by Lemma \ref{alpha-om} they bound a rectangle $Q$ in $S\sm A_j$.
The conjugate  harmonic functions  $\om_j$ and $\om_j^*$ are well defined on $Q$,
and $\om_j$ is constant on its horizontal sides, while $\om_j^*$ is constant on the vertical sides. 
Hence $\om_j + i\om_j^*$ is a conformal map of $Q$ onto a standard rectangle,
so that neither $\om_j$ nor $\om_j^*$ has critical points in $Q$. 
It follows that $Q$ is contained in one of the rectangles $\Pi_j^m$.

 A {\it harmonic rectangle} in $S$ is a subrectangle of some $\Pi_j^m$
saturated by the leaves of $\FF_j$. 

\begin{figure}[htbp]
\begin{center}
\makeabbrevs
\input{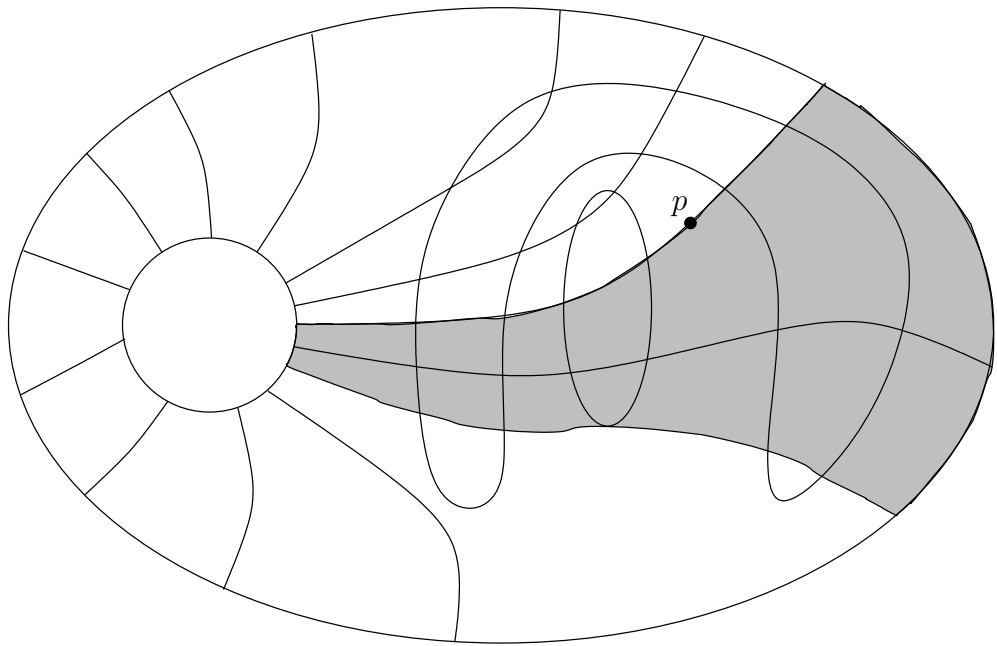}
\caption{Harmonic foliation $\FF_i$. Here $\BS$ and all $A_j$ are disks.
  The artificial puncture $p$ is made in  $\BS$ to  ensure that $S\sm A_i$ is not an annulus.
  One harmonic rectangle is shaded.}
\label{vertical_foliation}
\end{center}
\end{figure}

Any non-singular leaf $\CC$ of a harmonic foliation $\FF_j$ represents a good path in $S$.
Notice that the route $\RR(\CC)$ determines
the properly homotopy class of $\CC$ in $S\sm A_j$,
and hence determines the foliation $\FF_j$ and  the  rectangle $\Pi_j^m$ containing $\CC$.
These remarks, together with   Lemma~\ref{good rectangles} imply that
 the leaves with the same route,  $\RR(\CC)= \ba$, 
form  a (non-closed) harmonic rectangle  $P(\ba)$ in $S$.  
By Lemma \ref{route count},
there are at most $s(\Top, N)$ such routes $\ba$.
Therefore there are at most $Ns$ routes for the harmonic foliations 
 to all of the $N$ archipelagoes. 

 Associated big and little rectangles, $P_j(\ba)$ and $Q_j(\ba)$, $j=1,\dots l$, 
come together with any harmonic rectangle. 

\subsection{Buffers and the Small-Overlapping Principle}\label{so principle}

We are going to make  use of an important principle saying that 
{\it two wide path families have a relatively small overlap}.

A path family $\La$ on a rectangle $P$ is called a {\it genuinely vertical lamination} 
if the paths of $\La$ are genuinely vertical in  $R$, 
and the union of these paths, $\supp\La$, is measurable.
The projection to the horizontal side of $P$ (after uniformization by a standard rectangle)
induces a transverse measure $\nu$ on $\La$ (defined up to scaling).
If $P$ is embedded into a Riemann surface $S$  and $\gamma$ is a path on $S$, 
we say that  $\gamma$  {\it intersects less than $\eps$-portion of the total width of $\La$} if
$$
   \nu \{ \la\in\La :\; \la\cap \gamma \not=\emptyset\} < \eps \nu (\La)
$$
(note that this condition does not depend on the normalization of $\nu$).
The same discussion applies to the case of annulus. 

\begin{lem}\label{small overlap}
  Let $\kappa\geq 1$. 
 Let us consider a genuinely vertical lamination $\La$  on some annulus or rectangle $R\subset S$,
and let  $\GG$ be another path family  on $S$. If $\WW(\La)> \kappa$ and $\WW(\GG) \geq  \kappa$,
then there exists a path $\gamma\in \GG$ that intersects less than $1/\kappa$-portion of the
total width of $\La$.
In particular, if $\kappa=1$ then there is a path $\gamma\in \GG$ that does not cross
some leaf of $\La$.
\end{lem}

\begin{pf}
  Assume for definiteness that $R$ is a rectangle.
Let $\phi: E \ra R$ be the uniformization of $R$ 
 by a standard  rectangle $E = [0,a]\times [0,h]$ normalized so that
the projection of $\phi^* \La$ (which is a genuinely vertical lamination in $E$)
onto $[0,a]$ has length $\kappa$. 
Let us use the Euclidean metric $\mu$  on $E$ to bound $\WW(\La)$:
$$
    \WW(\La)\leq \frac{\area (\phi^* \La)} {\mu(\phi^* \La)^2}=\frac{\kappa}{h}
$$
(where $\area (\phi^* \La)$ stands for the area of $\supp \phi^*\La$). 
Since $\WW(\La) > \kappa$,
we conclude that $h < 1$, and thus $\area(\phi^*\La) <  \kappa$.

To bound $\WW(\GG)$, let us  use the push-forward metric $\rho=\phi_*(\mu|\La)$ on $S$.
If a curve  $\gamma\in \GG$ 
intersects at least $1/\kappa$-portion of the total width of $\La$, then 
the projection of $\phi^{-1}(\gamma)\subset E$ to $[0,a]$ has length at least 1, and hence 
 $$
     \rho(\gamma)= \mu(\phi^{-1}(\gamma)) \geq 1.
$$
If this happened for every $\gamma\in \GG$ then we would have 
$$
          \WW(\GG) \leq \area_\rho (\La) = \area(\phi^* \La) < \kappa  ,
$$
  contradicting the assumption.   
\end{pf}

Take some number $M>8$. 
Given a harmonic rectangle $P(\ba)$ of width greater than $M$, let us define two {\it buffers},
 $B^l(\ba)\subset P(\ba)$ and $B^r(\ba)\subset P(\ba)$,
as harmonic rectangles of width $M/2$ attached to the vertical sides of $P(\ba)$.   

\begin{lem}\label{disjoint curves}
 Let us consider two harmonic rectangles $P(\ba)$ and $P(\bbe)$ of width greater than $M$.
Then one can select four disjoint vertical leaves,
 one from each of the corresponding four buffers.
\end{lem}

\begin{pf}
Let $\La$ be the vertical foliation in $B^l(\ba)\cup B^r(\ba)$,
 and let  $\SS$ be the vertical foliation of $B_l(\bbe)$.  
  Applying  the previous lemma to  this data, we conclude that
there is a vertical leaf $\Gamma^l(\bbe)$ in $\SS$ that crosses less than  $1/4$ of the total width of $\La$.
Hence it crosses less than  1/2 of the total width of each $B^l(\ba)$ and $ B^r(\ba)$. 

Similarly, there is a vertical leaf $\Gamma^r(\bbe)$ that crosses less than
1/2 of the total width of each $B^l(\ba)$ and $B^r(\ba)$.
Together,  $\Gamma^l(\bbe)$ and  $\Gamma^r(\bbe)$ cross less than full width of each $B^l(\ba)$ and $B^r(\ba)$.
Hence  each $B^l(\ba)$ and $B^r(\ba)$ contains a vertical leaf, $\Gamma^l(\ba)$ and $\Gamma^r(\ba)$ respectively,
disjoint from both  $\Gamma^l(\bbe)$ and  $\Gamma^r(\bbe)$. 
\end{pf}

\subsection{Truncated rectangles and Disjointness Property}\label{truncated rectangles}

Let us remove the buffers from our harmonic rectangles:  
$$ 
     \tl P(\ba)= \cl (P(\ba) \sm (B^l(\ba)\cup B^r(\ba))).
$$
The associated truncated big and little rectangles will be naturally marked with tilde:
$\tl P_j(\ba)$ and  $\tl Q_j(\ba)$.

We can now formulate the key disjointness property for the truncated rectangles:

\begin{lem}\label{disjointness-2}
  If two associated truncated little rectangles $\tl Q_j(\ba)$ and $\tl Q_k(\bbe)$ overlap
then they represent the same proper homotopy class in $S\sm \La$ (up to orientation).
If their orientations match,
then one of the routes,  $\ba$ or  $\bbe$, is an extension of the other,
and $j=k$.
\end{lem} 

\begin{pf}
 Let us select in the buffers of  $P_j(\ba)$ and $P_k(\bbe)$
two disjoint pairs of leaves (by Lemma \ref{disjoint curves})
and consider the rectangles 
${\bf P}_j(\ba)\subset P_j(\ba)$ and ${\bf P}_k(\bbe)\subset P_k(\bbe)$
bounded by the corresponding pairs. 
By Lemma~\ref{inclusions},  their associated little rectangles, ${\bf Q}_j(\ba)$ and ${\bf Q}_k(\bbe)$,
contain the respective little rectangles  $\tl Q_j(\ba)$ and $\tl Q_k(\bbe)$.
Hence ${\bf Q}_j(\ba)$ and ${\bf Q}_k(\bbe)$ overlap as well. 
Since the big rectangles ${\bf P}_j(\ba)$ and ${\bf P}_k(\bbe)$ have disjoint vertical boundaries,
we can apply Lemma~\ref{disjointness} and 
complete the proof. 
\end{pf}


\begin{cor}\label{same itinerary}
   For any route $\ba$, the little rectangles $\tl Q_i(\ba)$ are pairwise disjoint.
\end{cor}

\begin{pf}
   Assume $Q_i(\ba)\cap Q_j(\ba)\not=\emptyset$ for some $i<j$.
Then by the first assertion of the lemma, one component of $\di Q_i(\ba)$ would lie on
$\di A_j$, which is impossible.  
\end{pf}

\begin{cor} \label{same height}
Suppose that $\tl Q_j(\ba)$ and $\tl Q_k(\bbe)$ overlap with matched vertical orientation.
Then $j = k$; 
moreover, if $|\ba|=|\bbe|$,
then $\ba = \bbe$.
\end{cor}

Fix your favorite $\de\in (0,1)$, e.g., $\de=1- \sqrt{2/3}$.
The total width of the rectangles $P(\ba)$ is equal to the modulus $Y$  
(by definition (\ref{Y}),  Example \ref{extremal foliation} and the Parallel Law).
For every route $\ba$,
we find that $\WW(\tl P(\ba)) \ge \WW(P(\ba)) - M$.
The number of routes $\ba$ is bounded by $Ns=Ns(\Top, N)$.
Therefore,
if $Y >  M Ns/\de $ then
\begin{equation}\label{tilde P}
     \sum_\sba \WW(\tl P (\ba)) >  (1-\de) Y.
\end{equation}

\subsection{$a$- and $b$-moduli}\label{forest}

We let 
$$
a_k = \sum_{|\sba|=k} \WW(\tl P(\ba)),
$$
$$
b_i^k = \sum_{|\sba|=k} \WW(\tl Q_i(\ba)),
$$
and $b_i = \max_{k \ge i} b_i^k$, $a = \sum_k a_k$, and $b = \sum_i b_i$.

We let 
$$ 
        x\oplus y=\frac{1}{\frac{1}{x}+ \frac{1}{y}} 
$$
be the {\it harmonic sum} of two numbers.

\begin{lem}\label{oplus}
The $a$- and $b$-moduli are related by the Series Inequality:
$$ 
       a_k \leq  \bigoplus_{i=1}^k  b_i.
$$  
\end{lem}
\begin{pf}
By Lemma \ref{good rectangles}, 
for each $\ba$ with $|\ba| = k$,
every vertical path of $\tl P(\ba)$ overflows each of  the little rectangles $\tl Q_i(\ba)$, 
with $1\le i \le k$. 
Moreover, by Corollary \ref{same itinerary}, the $\tl Q_i(\ba)$ are disjoint. 
Therefore,
by Lemma \ref{massive parallel},
$$
\sum_{|\sba|= k} \WW(\tl P(\ba)) \le \bigoplus_{i=1}^k \sum_{|\sba|=k} \WW(\tl Q_i(\ba)),
$$
and the Lemma follows.
\end{pf}

Let us now relate the $a$- and  $b$-moduli to the geometric moduli $X$,$Y$ and $Z$
in the Quasi-Additivity Law (see the Introduction). By  (\ref{tilde P}), 

\begin{equation}\label{a-moduli}
     a \geq (1-\de) Y,
\end{equation} 
provided $Y>MNs/\de$.
Furthermore,

\begin{lem}\label{b1-estimate}
    $ b_1  \leq X$. 
\end{lem}
\begin{pf}
We need to show that $b_1^k \leq X$ for every $k$.
Let us therefore fix $k$. 
By Corollary \ref{same height},
the $\tl Q_1(\ba)$ for $|\ba| = k$ are all disjoint,
so that the union of the associated vertical path families has width
 equal to 
$$
    \sum_{|\sba| = k} \WW(Q_1(\ba)) = b_1^k.
$$
On the other hand, 
this union is a subfamily of the family of paths connecting 
$\partial S$ and $\partial \Lambda$ in $S \sm \Lambda$ (recall that $\La=\bigcup A_j$).
Therefore 
$$
    \sum_{|\sba| = k} \WW(Q_1(\ba)) \le \WW(S, \Lambda) = X.
$$
\end{pf}

\begin{lem}\label{b-estimate}
    $  b  \leq Z. $ 
\end{lem}
\begin{pf}
  Let us arbitrarily label the archipelagoes $\{A_1, \ldots, A_n\}$.
  Let $\ba[i]$ denote the label of the $i^{\text{th}}$ archipelago 
  visited on the route $\ba$.
  We let 
$$
b_i^k(l) = \sum_{|\sba| = k;\, \sba[i] = l} \WW(\tl Q_i(\ba)),
$$
so that $b_i^k = \sum_l b_i^k(l)$.
Let $k : \N \to \N$ be such that $b_i = b_i^{k(i)}$.

We claim that 
$$
\sum_i b_i^{k(i)}(l) \le \WW(S \sm \bigcup_{k\neq l} A_k, A_l);
$$
this would imply (by summing over $l$) that $b \le Z$.
To show the claim, 
first note that the $\tl Q_i(\ba)$ for $|\ba| = k(i)$ and $\ba[i] = l$
are disjoint (where $l$ is fixed and $i$ is arbitrary):
any two such rectangles have the same roof,
so they have the same vertical orientation if they overlap; 
then by Corollary \ref{same height}, 
they have the same height $i$ and therefore the same route $\ba$.
Moreover the vertical paths of these $\tl Q_i(\ba)$ all connect 
$\partial(S \sm \bigcup_{k\neq l} A_k)$ to $\partial A_l$ in $S \sm \Lambda$;
the claim follows.
\end{pf}

\subsection{An arithmetic inequality}\label{arithmetic sec}

\begin{lem}\label{harmonic sums} 
   Consider two sequences of positive numbers, $\{a_i\}_{i=1}^n$ and $\{b_i\}_{i=1}^n$,
 such that $a_1= b_1$, $a_{i+1}\leq  a_i\oplus b_{i+1}$. Then
\begin{equation}\label{inequality}
         \left(\sum_{i=1}^n  a_i\right)^2 \leq \frac{4}{3} \, b_1 \sum_{i=1}^n b_i. 
\end{equation}
\end{lem}
         
\begin{pf} 
Let 
$$
a = \sum_{i=1}^n a_i, \quad b = \sum_{i=1}^n b_i.
$$
We have, for $i > 1$,
$$
b_i \ge \frac{a_ia_{i-1}}{a_{i-1}-a_i} = a_i + \frac{a_i^2}{a_{i-1} - a_i}.
$$
and therefore
\begin{eqnarray}
b - b_1 &\ge& \sum_{i = 2}^n  \left( a_i + \frac{a_i^2}{a_{i-1} - a_i} \right) \nonumber\\
&=&
a - a_1 + \sum_{i=2}^n \frac{a_i^2}{a_{i-1} - a_i} \nonumber\\
&\ge& \label{Cauchy-Schwarz}
a - a_1 + \frac{\big(\sum_{i=2}^n a_i \big)^2}{\sum_{i=2}^n (a_{i-1} - a_i) } \\ 
&=&
a - a_1 + \frac{(a-a_1)^2}{a_1 - a_n} \nonumber \\
&\ge&
a - a_1 + \frac{(a - a_1)^2}{a_1},  \nonumber
\end{eqnarray}
where inequality (\ref{Cauchy-Schwarz}) follows from the Cauchy-Schwarz inequality written as follows:
$$
\big(\sum x_j\big)^2 \le \sum \frac{x_j^2}{y_j} \sum y_j.
$$
Therefore, because $a_1 = b_1$,
$$
\frac{b_1 b}{a^2} \ge 1 - \frac {a_1} {a} + \left( \frac {a_1} {a} \right)^2 \ge \frac34.
$$

\end{pf}

  \subsection{Completion of the proof of the QA Law}
Let us consider the $a$- and $b$-moduli from  \S \ref{forest}. 
Lemma  \ref{oplus} puts us into a position to apply  estimate (\ref{inequality})
to these moduli.
Incorporating (\ref{a-moduli}) and Lemmas  \ref{b1-estimate} and \ref{b-estimate} into (\ref{inequality}),
we obtain:
$$
   (1-\de)^2 Y^2  \leq     \frac{4}{3} XZ,    
$$
provided $Y> Ms/\de$, and we are done.
\QED

\subsection{QA Law: Variations}\label{variation}
We will now formulate several variations and special cases of the QA Law
suitable for the dynamical applications. 

\subsubsection { Fractal archipelagos.}  \label{fractal}
A compact set $A \subset \inter S$ is called a set of {\it finite type},
or a {\it (closed) archipelagos}, if $A = \cap U_i$
where $U_i$ is a nested sequence of open archipelagos of bounded topological complexity
(equivalently, $S\sm K$ is a Riemann surface of finite type).
In this case, we let
$$
    \Top_S(A)= \liminf \Top_S(U_i).
$$  
(If we have a finite family of disjoint closed archipelagos $A_j$,
we let $\Top_S \{A_j\} = \Top_S(\bigcup A_j)$.) 

{\it By an approximation argument, the QA Law is valid for
these more general archipelagos as well.}

\subsubsection { Collars.}\label{QA Law with collars} 
Let  $A_j' $ be a topological disk such that 
$$
   A_j\subset A_j' \subset S\sm \bigcup_{k\not=j} A_k.
$$
If $\mod (A_j', A_j)\geq \eta\mod(S, A_j) >0$, then we call $A_j'$ an {\it $\eta$-collar} around $A_j$.    
If all the archipelagos $A_j$ have $\eta$-collars, 
we say that the archipelagos satisfy the {\it $\eta$-Collar Assumption}.
Under this assumption, they are $\eta^{-1}$-separated
(since $Z\leq \sum \WW(A_j', A_j)$). Thus, we obtain:  

\proclaim QA Law with collars.
Under the  $\eta$-Collar Assumption,
there exists $K$ depending only on $\eta$ and $\Top_S\{A_j\}$
 such that:
 $$ Y\geq K \imply Y\leq 2\eta^{-1} X. $$ 

One can also allow general {\it holomorphic collars} instead of embedded ones.
Precisely speaking, assume $A_j$ is embedded into an abstract conformal disk $A_j'$ which in turn is 
mapped into $S\sm \bigcup_{k\not= j} A_k $  holomorphically by some  map $i$
such that $i| A_j=\id$ and $i^{-1}(A_j)=A_j$.  
If $\mod (A_j', A_j)\geq \eta \mod (S, A_j) >0$, then we call $A_j'$ a {\it holomorphic} $\eta$-collar around $A_j$.
Since every path connecting $A_j$ to the rest of the boundary of $S\sm \bigcup A_k$
can be lifted to a vertical path in $A_j'\sm A_j$,  
Corollary \ref{two families} yields: $Z\leq \sum \WW(A_j', A_j)$. 
Thus, the {\it $\eta$-Collar Assumption for holomorphic collars
implies $\eta^{-1}$-separation of the archipelagos as well}. 

\subsubsection{Comparable terms}

In further applications in  holomorphic dynamics,
we will often encounter the situation when the individual terms that
appear in  the moduli  $Y$ and $Z$ are all comparable.
Here is the user-friendly version of the Quasi-Additivity Law
in this situation:

\proclaim  QA  Law with comparable terms.
    Fix some $\eta\in (0,1)$. 
Let $W\Subset \inter U$ and  $D_i'\Subset \inter W $, $i=1,\dots, N$, 
be topological disks  such that the closures of $D_i'$ are pairwise disjoint,
and let $D_i\Subset D_i'$ be smaller disks. 
Then there exists a $\de_0 > 0$ (depending on $\eta$ and $N$) such that: \\
If  for some $\de\in (0, \de_0)$ and for all $i$, 
$$
           \eta\de <  \mod(D_i'\sm D_i)  \leq \mod (U \sm D_i) < \de,   
$$
    then
$$
  \mod(U \sm W) < \frac{2\eta^{-1} \de}{N}.
$$

Of course, this version is a particular case of the QA Law with collars.

\section{Quasi-Invariance Law}

In this section, we
will prove a general transformation law for conformal moduli under covering maps.
To this end, we will make use of the following well-known result:

\begin{prop} \label{galois}
Let $f: U \to V$ be a branched cover of Riemann surfaces of degree $N$. 
Then there is a Galois branched cover $g : S \to V$ of degree at most $N!$ 
that factors as $g = f \circ h$ for some $h: S \to U$.
Moreover, $g$ is ramified only over critical values of $f$. 
\end{prop}

The proof uses a lemma that is a simple exercise in group theory:

\begin{lem}\label{group action}
Suppose that $H$ is a subgroup of a group $G$, and \break $[G:H] = N$. 
Then there is a normal subgroup $L$ of $G$ such that $L < H$, and $[G:L] \le N!$.
\end{lem}

\begin{pf}
The coset action of $G$ on $G/H$ provides a homomorphism from $G$ to the group of permutations of $G/H$, 
which has order at most $N!$. We let $L$ be the kernel of this homomorphism; it has the desired properties.
\end{pf}

\nin {\it Proof of Proposition \ref{galois}}. 
Let $O$ be the set of critical values of $f$, and let $E=f^{-1}(O)$.
Then $f: U \sm E \to V \sm O$ is an unbranched cover of degree $N$. 
Hence  $f_*\pi_1(U\sm E)$ has index $N$ in $\pi_1(V\sm O)$, 
so by Lemma~\ref{group action} we can find a subgroup of $f_*\pi_1(U\sm E)$ that is a  normal subgroup of $\pi_1(V\sm O)$
of degree at most $N!$. There is then the corresponding cover $g: S' \to V\sm O$ which we can complete to a branched cover 
$g: S \to V$ with the desired properties. 
\QED

 We say that a closed set $K\subset S$ is a {\it hull}
if it is a full connected non-degenerate continuum.

Given a holomorphic map $f: S\ra S'$, and two closed subsets $K\subset S$, $K'\subset S'$ such that
$f(K)\subset K'$, we say that the restriction $f: K\ra K'$ is a {\it branched covering of degree $d$} if:

\ssk \nin $\bullet$ 
   For any $x\in K$, there exists a  neighborhood $U\ni x$ such that $K\cap U = f^{-1}(K')\cap U$;

\ssk \nin $\bullet$ 
  For  any point regular value $x'\in K'$ of $f$, $\# (f|K)^{-1}(x)=d$.  
 
\ssk
Let us consider a Riemann surface $V$ with several archipelagos $B_j$ contained in hulls $B_j'$,
and several marked points $v_i$ (some of them may belong to the archipelagos or the hulls).  
Let us consider the family $\GG_k$ of proper curves 
$\gamma\subset S\sm \bigcup B_j $ that begin on $B_k$ and satisfy one of the following conditions:

\ssk \nin $\bullet$
   $\gamma$ ends on another archipelagos  $B_j$, $j\not=k$;

\ssk \nin $\bullet$
 if $\gamma$ ends on the same $B_k$ then it does not pass through the marked points $v_i$ 
  and is  {\it non-trivial} in the sense that it cannot be homotopic in  $S\sm (B_k\cup \{v_i\})$
  to an arbitrary small neighborhood  of the hull $B_k'$.%
\footnote{Notice that such a trivial $\gamma$ is allowed to have arbitrary complexity in $S\sm B_k$.}  

\ssk \nin
Under these circumstances, we let  
$$
     Z_S\{ B_j , B_j',  v_i\}  = \sum_k \WW(\GG_k). 
$$

\proclaim General Quasi-Invariance Law.
Let us consider the following data:
\begin{itemize}
\item         Two Riemann surfaces of finite type, $U$ and $V$;
\item 
     Two closed sets $\displaystyle{\La'=\bigcup_{j=1}^p \La_j'\subset U}$ and 
$\displaystyle{B'=\bigcup_{j=1}^p B_j' \subset V} $ whose connected components, 
$\La_j'$ and $B_j'$ respectively, are hulls;
\item
   Two families of compact archipelagos,  $\La_j\subset \La_j'$ and $B_j\subset B_j'$;
\item
     A branched covering $f: U \ra V$ of degree $D$ that restricts to branched coverings
$f: \La_j'\ra B_j'$ of degree $d_j$. Suppose $\La_j$ is the union of some components of $f^{-1}(B_j)$,
and let $\Crit$ stand for the set of critical values of $f$.
\end{itemize}
Then there exists $K$ depending on $\Top_V\{B_j\}$ and $D$ such that
     $$ Y_U\{\La_j\} >K \imply  Y_U \{\La_j\}^2 \leq 2 d^2\, X_V \{B_j\} \, Z_V \{B_j , B_j' , \Crit\}. $$

\begin{pf}
If we replace the archipelagos $\La_j$ with  $\La_j = (g|\La_j')^{-1}(B_j)$
we make the left-hand side bigger without changing the right-hand side.
So, we can assume without loss of generality that $\La_j = (g|\La_j')^{-1}(B_j)$

Let $\Crit \subset V$ be the set of critical values of $f$,
and let $E=f^{-1}(\Crit)\subset U$.
By Proposition \ref{galois}, there exists a  branched covering $h: S\ra U$ of degree at most $(D-1)!$               
 with critical values in $E$
such that $g = f\circ h: S\ra U$ is a Galois branched covering.
Let $\Gamma$ be the Galois group of the covering $g$ acting on $S$. 

Let $A_j'(i) \subset S$ be the connected components of $g^{-1} (B_j')$
labeled in such a way that $h(A_j'(1))= B_j'$, and let $ A_j' =  A_j'(1)$. 
For any given $j$, these components are transitively permuted by $\Gamma$. 
We let $L_j$ be the number of these components.
 
Let us also consider the corresponding archipelagos 
$$
   A_j(i)  = (g|\, A_j'(i) )^{-1} (B_j), \quad A_j \equiv A_j(1), 
$$
$$  
                A=\bigcup A_j(i) = g^{-1} ( B).
$$

Let $\XX$, $\YY$ and $\ZZ$ stand respectively for the $X$-, $Y$- and $Z$-moduli for this family of archipelagos.
By Lemma \ref{modulus transform-1} from the Appendix, we have:
\begin{equation}\label{XXX}
   \XX= |\Gamma|\,  X_V\{\La_j\}.
\end{equation}

Let $m_j=\deg (h: A_j' \ra  \La_j' )\equiv \deg (h: A_j\ra \La_j) $. 
Then the stabilizer of $A_j'$ in $\Gamma$ consists of $d_j m_j$ elements,
and hence the $\Gamma$-orbit of $A_j$ consists of $L_j=|\Gamma|/d_j m_j$ archipelagos $A_j(i)$. 
Since for each $j$, these archipelagos are symmetric in $S$, we have: 
\begin{equation}\label{YYY}
  \YY= \sum_j \frac {|\Gamma|}{d_j m_j} \WW(S, A_j) \geq \frac {|\Gamma|}{d}  \sum_j \WW(U, \La_j)= 
                     \frac{|\Gamma|}{d} Y_U\{B_j\}
\end{equation}
where the middle inequality follows from Lemma \ref{modulus transform-2}.


\msk Let us now show that
\begin{equation}\label{ZZZ}
   \ZZ \leq { |\Gamma|\,  Z_V\{A_j, \Crit\}}+C,
\end{equation}
where $C$ depends only on $(\Top_U\{A_j\})$ (which in turn depends only on $\Top_V\{B_j\}$ and $D$). 

For any $k\in [1,p]$, let us consider the harmonic foliation $\FF_k$ that measures the extremal
width between $A_k$ and the rest of the boundary of $S\sm \bigcup A_j(i)$ (see \S \ref{Vert fol}).
Then $S\sm \bigcup A_j$ is tiled by the harmonic rectangles $\Pi^n_k$, $n=1,\dots, s_k$. 
Their total  number $\sum s_k$ depends only on $\Top_S\{A_j\}$.  
Applying group $\Gamma$, we obtain a family of harmonic  rectangles
$\Pi^n_{ji}$ (connecting $A_j(i)$ to the rest of the boundary of $S\sm \bigcup A_j(i)$)
that are permuted by the $\Gamma$-action.

Let $\tl \Pi^n_{ji}$ be the truncated rectangle obtained
by removing two buffers of width four each
from $\Pi^n_{ji}$ (as in \S \ref{truncated rectangles}).
They are also permuted by $\Gamma$.
Since these rectangles represent different homotopy classes in $S\sm \bigcup A_j(i)$,
Lemma \ref{disjoint curves} implies that
the truncated rectangles are pairwise disjoint.

Since the fibers of $g$ coincide with the orbits of $\Gamma$,
each $\tl \Pi_k^n$ projects injectively onto some proper rectangle $\tl Q_k^l$ in $V\sm \bigcup B_j$,
and these rectangles are pairwise disjoint. Moreover, there are $d_k m_k$ rectangles $\Pi_k^n$
that project onto each $\tl Q_k^l$. 

The foliation $\FF_k$ on $\bigcup_n \tl\Pi^n_k$ descends to a foliation
$\HH_k$ supported on $\bigcup_l \tl Q^l_k$. The leaves of this foliation belong to the 
family $\GG_k$ of curves defining the modulus $Z_V\{B_j, \Crit\}$. 
(Indeed, if some leaf $\gamma$ connecting $B_k$  to itself was trivial then
it would lift to a path connecting  $A_k'$ to itself). Hence
$$
     \sum_l \WW(\tl Q^l_k) \leq \WW(\GG_k),
$$    
  and we obtain:
$$
     \WW(S\sm \bigcup_{j\not=k}A_j, \, A_k) =\sum_n \WW(\tl \Pi^n_k) + 8 s_k 
$$
$$
\leq d_k m_k\sum_l \WW(\tl Q^l_k) + 8 s_k \leq  d_k m_k \WW(\GG_k)  + 8 s_k .
$$
(Here $8 s_k$ appears as the total width of the buffers removed.)
Multiplying  the last estimate by $L_k$ and summing  up over $k$
(making use of the symmetry and of $|\Gamma|=L_k d_k m_k$), 
we obtain (\ref{ZZZ}).

\msk
By the Quasi-Additivity Law, $\YY^2\leq $2.75$ \XX \ZZ$.
Together with (\ref{XXX}) and (\ref{YYY}) and (\ref{ZZZ}) it implies the desired estimate,
provided $\ZZ$ is sufficiently big (which is certainly the case when $Y_V\{\La_j\}$ is 
sufficiently big).
\end{pf}

\subsection{QI Law: Variations}\label{qi variations}
  Let us now list several variation and special cases of the General QI Law. 
In what follows, 
the setting of the General QI Law is assumed,  
and we let $X_U= X_U\{\La_j\}$, $Y_V=Y_V\{B_j\}$, $Z_V= Z_V\{B_j, B_j', \Crit\}$.  

\subsubsection{Quasi-Additivity Law} It is a particular case of the QI Law when $f=\id$.

\sss{Separation}
In the context of the QI Law, the $\xi$-Separation Assumption should be formulated as follows:
$$
   Z_V  \leq \xi Y_U .
$$

\proclaim QI Law with separation.
If  the archipelagos $B_j$  are $\xi$-separated,
then there exists $K$ depending only on $\xi$, $\Top_V\{B_j\}$ and $D$ 
 such that:
 $$ Y_U  \geq K \imply  Y_U \leq 2 \xi d^2 X_V. $$

\sss{Collars}
  The definition of  $\eta$-collars should  also be adjusted in this more general context.
Namely, a disk  $\BB_j \supset B_j'$ is called an {\it $\eta$-collar } of
$B_j$ if $\BB_j \sm B_j'\subset V\sm (\bigcup_{k\not=j} B_k\cup \Crit )$ and
\begin{equation}\label{collar assumption}
   \mod (\BB_j , B_j)\geq \eta \mod (U, \La_j).
\end{equation}
More generally, one can define a {\it holomorphic $\eta$-collar} $\BB_j$  as an abstract conformal disk
$\BB_j$ such that $B_j'$ is embedded into $\BB_j$ and there is a holomorphic map
$i: \BB_j\ra V$ such that $i|B_j'=\id$,   
$$
   i(\BB_j\sm B_j')\subset V\sm (\bigcup B_k\cup \Crit), 
$$ 
and (\ref{collar assumption}) is satisfied.

\proclaim QI Law with collars.
If all the the archipelagos $B_j$ have holomorphic $\eta$-collars then
there exists $K$ depending only on $\eta$,  $\Top_V\{B_j\}$ and $D$
 such that:
 $$ Y_U \geq K \imply Y_U  \leq 2 \eta^{-1} d^2 X_V . $$ 

\sss{No critical values outside $\bigcup B_j$}

In what follows  we will restrict ourselves to the case when $U$ and $V$ are topological disks,
and each archipelago consists of  a single island.

If $\Crit\subset \bigcup B_j$ then the family $\GG_k$  of curves defining $Z_V\{B_j, \Crit \}$
is contained in the family $\hat\GG_k$ of all non-trivial proper curves $\gamma$ in $V\sm \bigcup B_j$
that begin on $B_k$ (so $\gamma$ is  allowed to end on $B_k$ as long as it goes around some
other island). Let 
$$  
   \hat Z_V = \hat Z_V\{B_k\}=\sum \WW(\hat \GG_k\}.
$$
Then  $ Z_V \leq \hat Z_V $,  so that the condition 
\begin{equation}\label{hat-separation}
  \hat Z_V \leq \xi Y_U .
\end{equation}
is stronger than $\xi$-separation.
This leads us to the following variation of QI Law:

\proclaim QI Law with all critical values in $\bigcup B_j$. 
For any natural numbers $p, D\geq d$ and any  $\xi>0$, 
there exists a $K=K(p,D,\xi)$ with the following property.  
Let $U$ and $V$ be two topological disks,
and let $\La_j$ and $B_j$ be two families of disjoint compact hulls in $\inter U$ and $\inter V$ 
respectively, $j=1,\dots, p$.
Let  $f\from (U,  \{\La_j\} )\to (V, \{ B_j\})$  be a branched covering 
with critical values in $\bigcup B_j$ such that 
 $\deg(f\from U\to V) = D$,
$\La_j$ is a component of $f^{-1}(B_j)$, and  $\deg (f\from \La_j\to B_j)\leq d$, $j=1,\dots, p$.
Under the separation assumption (\ref{hat-separation}) 
we have: 
$$
    Y_U  > K \imply Y_U   \leq 2\xi d^2 X_V .
$$

\sss {Covering Lemma.} The Basic Covering Lemma
  stated in the Introduction is a special case of the General QI Law with embedded collars
 when both Riemann surfaces, $U$ and $V$, are conformal disks, and the archipelagos $\La$ and $B$ consist
of a single Jordan island each. 
In the following variation the collars are allowed to be holomorphic:  

 \proclaim Covering Lemma with holomorphic collars.
Fix some  $\eta\in (0,1)$. Let us consider
two topological disks $U$ and $V$, two hulls $\La'\subset U$ and $B'\subset V$,
and two compact hulls $\La\subset \La'$ and $B\subset B'$.
\\
  Let $f: U \ra  V$ be a branched covering of degree $D$ such that
$A'$ is a component of $f^{-1}(B')$, and $\La$ is the union of some components of $f^{-1}(B)$.
Let $d=\deg(f : \La'\ra B')$.\\
Assume  $B'$ is also embedded into a holomorphic $\eta$-collar $\BB'$, i.e., 
 there is a holomorphic map $i: \BB \ra V$ such that
$i|\, B'=\id$, $i^{-1}(B')=B'$,
 and $i(\BB) \sm B'$ does not contain the critical values of $f$  
such that
$$
   \mod(\BB, B) > \eta \mod(U,  \La).
$$
Then 
$$
   \mod(U, \La) < \eps(\eta, D)  \imply   \mod (V, B)  <  2\eta^{-1} d^2 \mod(U,  \La).
$$

The Basic Covering Lemma stated in the Introduction is used in \cite{high},
the Covering Lemma with holomorphic collars  is used in \cite{decorations}, while
the QI Law with all critical values in $\bigcup B_j$  is used in  \cite{K}.

\section{Appendix: Extremal length and width}

There is a wealth of sources containing background material on extremal length,
see, e.g., the book of Ahlfors \cite{A}. 
We will briefly summarize the necessary minimum. 

\subsection{Definitions}
Let $\GG$ be a family of curves on a Riemann surface $U$.
Given a (measurable) conformal metric $\mu= \mu(z) |dz|$ on $U$, 
let
$$   \mu (\GG) = \inf_{\gamma\in \GG}  \mu (\gamma).$$
 where $\mu(\gamma) $ stands for the $\mu$-length of $\gamma$.
The length of $\GG$  with respect to $\mu$ is defined as   
$$
   \LL_\mu (\GG) = \frac {\mu(\GG)^2 } {\mu^2 (U)}, 
$$
where $\mu^2= \mu(z)^2 dz\wedge \bar dz$ is an area form of $\mu$.
Taking the supremum over all conformal metrics $\mu$, we obtain the {\it extremal length}
$\LL(\GG)$ of the family $\GG$. 

  The {\it extremal width} is the inverse of the extremal length:
$$
 \WW(\GG)=\LL^{-1}(\GG).
$$
It can be also defined as follows. Consider all conformal  metrics $\mu$ such that
$\mu(\gamma)\geq 1$ for any $\gamma\in \GG$.  Then $\WW(\GG)$ is the infimum of 
the areas $\mu^2(U)$ of  all such metrics.

\begin{example}\label{extremal foliation}
    For a standard rectangle $P=I\times [0,h]$, let $\GG$ be the family of vertical curves,
and let $\La$ be the genuinely vertical foliation.
Then 
$$ 
    \LL(\GG)= \LL(\La) = \frac{h}{|I|} \equiv \mod P.
$$
The similar formulas hold for the standard cylinder $C = \T\times [0,h]$.  
\end{example}

\subsection{Electric circuits laws} 
We say that a family  $\GG$ of curves {\it overflows} a family $\HH$ if
any curve of $\GG$ contains some curve of $\HH$.
We say that two families, $\GG_1$ and $\GG_2$, are disjoint if
any  two curves, $\gamma_1\in \GG_1$ and $\gamma_2\in \GG_2$, are disjoint.

We let $ x\oplus y = (x^{-1} + y^{-1})^{-1}$ be the {\it harmonic sum} of $x$ and $y$
(it is conjugate to the usual sum by the inversion map $x\mapsto x^{-1}$). 

The following crucial properties of the extremal length and width show
that the former behaves like the resistance in electric circuits,
while the latter behaves like conductance.

\proclaim  Series Law/Gr\"otzsch Inequality. 
 Let $\GG_1$ and $\GG_2$ be two disjoint families of  curves,
 and let $\GG$ be   a third family
that overflows both $\GG_1$ and $\GG_2$. Then 
$$  \LL(\GG) \geq \LL(\GG_1) + \LL(\GG_2), $$
or equivalently,
$$  \WW(\GG)\leq \WW(\GG_1)\oplus \WW(\GG_2). $$ 

\proclaim Parallel Law.
  For any two families $\GG_1$ and $\GG_2$ of curves we have:
$$
    \WW(\GG_1\cup \GG_2)\leq \WW(\GG_1) + \WW(\GG_2).
$$
If $\GG_1$ and $\GG_2$ are disjoint then 
$$
  \WW(\GG_1\cup \GG_2) = \WW(\GG_1) + \WW(\GG_2)
$$

Note that the Parallel Law inequality implies the estimate $X\leq Y$ between the moduli
from the Introduction. 

\ssk

From the Series and Parallel Laws we can derive the following more general result:
{\def\l{\lambda}\def\L{\Lambda}\def\D{\Delta}\def\G{\Gamma}
\begin{prop} \label{massive parallel}
Suppose that $\D_\l^i$, $\G_\l$ for $i = i \ldots k$, $\l \in \L$
(where $\L$ is finite)
are path families supported on a Riemann surface $S$. 
Assume for each $\l \in \L$,
the $\D^i_\l$ have disjoint support,
and $\G_\l$ overflows each of the $\D^i_\l$.
Then
$$
\sum_\l \WW(\G_\l) \le \bigoplus_{i=1}^k \sum_\l \WW(\D_\l^i).
$$ 
\end{prop}
\begin{pf}
We form path families $\hat \D^i_\l$ and $\hat \G_\l$ 
 on the Riemann surface $S \times \L$ 
by putting $\D^i_\l$ and $\G_\l$
on the copy of $S$ labelled by $\l$.
Let $\hat \G = \bigcup_\l \hat \G_\l$
and  $\hat \D^i = \bigcup_\l \hat \D_\l^i$.
By the Parallel Law,
$$
   \WW(\hat \G) = \sum_\l \WW(\hat \G_\l), \quad
   \WW(\hat \D^i) = \sum_\l \WW(\hat \D^i_\l).
$$
Moreover,   $\hat \G$ overflows each of the $\hat \D^i$,
and the $\hat \D^i$ are disjoint.
Therefore,  
by the Series Law,
$$
\sum_\l \WW(\hat \G_\l) \le \bigoplus_{i=1}^k \sum_\l \WW(\hat \D_\l^i),
$$
and the result follows.
\end{pf}
}

\subsection{Transformation rules}
 Both extremal length and extremal width  are conformal invariants. 
 More generally, we have:
\begin{lem}\label{increase}
   Let $f: U\ra V$ be a holomorphic map between two Riemann surfaces,
and let $\GG$ be a family of curves on $U$. Then 
$$
      \LL(f(\GG))\geq \LL(\GG).
$$
Moreover, if $f$ is at most $d-to-1$, then 
$$
   \LL (f(G)) \leq d\cdot  \LL(\GG).
$$ 
\end{lem}

\begin{pf}
  Let $\mu$ be a conformal metric on $U$. 
Let us push-forward the area form   $\mu^2$  by  $f$. 
We obtain the area form $\nu^2= f_*(\mu^2)$  of some conformal metric $\nu$ on $V$.   
Then $\nu^2 (V) = \mu^2 (U)$ and $f^* (\nu) \geq \mu$. It follows that
$$
    \LL_\mu(\GG) \leq  \LL_\nu(f(\GG)) \leq  \LL(f(\GG)) .
$$ 
Taking the supremum over $\mu$ completes the proof of the first assertion. 

\ssk 
  For the second assertion,
let us consider a conformal metric $\nu$ on $V$ and pull it back to $U$, $\mu= f^* \nu$.
Then $\mu(\gamma) = \nu(f(\gamma))$ for any $\gamma\in \GG$, while $\mu^2(U) \leq d\cdot \nu^2(V)$. 
Hence 
$$
   \LL(\GG)\geq \LL_\mu(\GG) \geq \frac {1}{d} \LL_\nu(f(\GG)),
$$
and taking the supremum over $\nu$ completes the proof. 
\end{pf}

\begin{cor}\label{two families}
     Under the circumstances of the previous lemma, 
let $\HH$ be a family of curves in $V$ satisfying the following lifting
property: any curve $\gamma\in \HH$ contains an arc that lifts to some curve
in $\GG$. Then $\LL(\HH)\geq \LL(\GG)$.
\end{cor}

\begin{pf}
   The lifting property means that the family $\HH$ 
overflows the family $f(\GG)$. Hence $\LL(\HH)\geq \LL(f(\GG))$,
and the conclusion follows.
\end{pf}

\subsection{Extremal distance and Dirichlet integral}

 Given a compact subset $K\subset \inter U$, the {\it extremal distance} 
               $$\LL(U, K)\equiv \mod (U,K)$$
(between $\di U$ and $K$) is defined as $\LL(\GG)$, where $\GG$ is the family of curves 
connecting $\di U$ and $K$. In case when $U$ is a topological disk and $K$ is connected,
we obtain the usual modulus $\mod(U\sm K)$ of the annulus $U\sm K$. 

\ssk {\it Remark.} $\LL(U,K)$ can also be defined as $\LL(\GG')$ where $\GG'$ is the family of curves 
{\it in} $U\sm K$ connecting $\di U$ to $K$. Indeed, since $\GG\supset \GG'$, $\LL(\GG)\leq \LL(\GG')$.
Since each curve of $\GG$ overflows some curve of $\GG'$, $\LL(\GG)\geq \LL(\GG')$. 
One can also make a compromise and use the intermediate  family of curves in $U$ connecting $\di U$ to $K$. 
\ssk

 We let $\WW(U,K) = \LL^{-1}(U,K)$.

\begin{lem}\label{modulus transform}
 Let $f: U\ra V$ be a branched covering of degree $N$ between two compact Riemann surfaces with boundary.
Let $A$ be a compact subset of $\inter U$ and let $B=f(A)$.  
Then 
$$
   \mod(U,  A)\leq \mod(V, B) \leq N \, \mod(U,A).
$$
\end{lem}

\begin{pf}
  Let $\GG$ be the family of curves in $U$ connecting $\di U$ to $A$,
and let $\HH$ be the similar family in $V$.  
Notice that every curve $\gamma\in \HH$ lifts to a curve in $\GG$:
begin the lifting on $A$, and it must end on $\di U$ since $f: U\ra V$ is proper. 
Thus,  $\HH= f(\GG)$, and  Lemma \ref{increase} completes the proof. 
\end{pf}

 Extremal width $\WW(U,A)$ can be explicitly expressed as the
{\it Dirichlet integral of the harmonic measure} (see \cite[\S 4-9]{A}):
$$
      \WW (U,A) = 4 \int_{U\sm A} |\di h|^2 
$$ 
where $h: V\sm B \ra \R$ is the harmonic function equal to 1 on $\di B$ and vanishing on $\di U$, 
and $|\di h|^2$ is the area form associated with the holomorphic differential $\di h = (\di h/ \di z) dz $.

\comm{**********************
\subsection{Dependence on the domain}

  Recall that the {\it Hausdorff distance} between two compact subsets $K, Q\subset S$
is defined as the infimum of numbers $\eps>0$ such that $K$ is contained in the $\eps$-neighborhood of $Q$
and the other way around.

\begin{lem}\label{continuity}
  Let $U_n$ be a sequence of domains in $S$ containing compact subsets  $K_n\Subset U_n$.
and let  $(U, K)$ be a similar data, where $U\Subset S$. 
If  $\di U_n\to \di U $ and $K_n\to K$ in the Hausdorff metric,
 then $\WW(U_n,K_n)\to \WW(U,K)$. 
\end{lem} 
************************************}

\subsection{More transformation rules} 

The Dirichlet integral formulation allows us to sharpen the lower bound in Lemma \ref{modulus transform}:

\begin{lem}\label{modulus transform-2}
 Let $f: U\ra V$ be a branched covering between two compact Riemann surfaces with boundary.
Let $A$ be an archipelago in $U$,  $B=f(A)$, and 
assume that  $f:  A\ra B$ is a branched covering of degree $d$. Then   
$$
      \mod(V, B) \geq d\, \mod (U,  A).
$$   
\end{lem}

\begin{pf}
The Riemann surface  $V \sm B$ is decomposed into finitely many rectangles saturated by the leaves of the harmonic flow
(see \S \ref{Vert fol}). Slit these rectangles by the leaves containing the critical values of $f$.
We obtain finitely many foliated rectangles $\Pi_i$ such that
$$
         \sum \WW(\Pi_i) = \WW(V, B).
$$

Each of these rectangles lifts to $d$ properly embedded rectangles $P_i^j$ in $U\sm A$
(with the horizontal sides on $\di U$ and $\di A$). Moreover, $\WW(P_i^j)=\WW(\Pi_i)$.
Hence
$$
        \WW(U,A) \geq \sum \WW(P_i^j) = d\,  \WW(V,B). 
$$
\end{pf}

{\it Remark.} A similar estimate is still valid for an arbitrary compact set $A$,
and can be proved by approximating $A$ by  archipelagos.

\ssk

Putting the above two lemmas together
(or using directly that the Dirichlet integral is transformed as the area     
under branched coverings) we obtain:  


\begin{lem}\label{modulus transform-1}
Let $(U,A)$ and $(V,B)$ be as above,
and let  $f: U\sm A\ra V\sm B$ be a branched covering of degree $N$. Then   
$$
      \mod(V, B)= N\, \mod (U,  A).
$$   
\end{lem}

\end{document}